\documentclass[10pt]{amsart}

\usepackage[ansinew]{inputenc}
\usepackage[T1]{fontenc}
\usepackage{amsmath}
\usepackage{amssymb}
\usepackage{oldgerm}
\usepackage[curve]{xypic}
\usepackage{amsthm}
\usepackage{hyperref}
\usepackage{stmaryrd}
\usepackage{graphicx}

\newtheorem{satz}{Satz}[section]
\newtheorem{Theorem}[satz]{Theorem}
\newtheorem{Lemma}[satz]{Lemma}
\newtheorem*{Hauptsatz}{Theorem}
\newtheorem{Prop}[satz]{Proposition}
\newtheorem{example}[satz]{Example}
\newtheorem{Cor}[satz]{Corollary}
\newtheorem{Conj}[satz]{Conjecture}
\newtheorem*{Conj*}{Conjecture}

\theoremstyle{definition}
\newtheorem{Definition}[satz]{Definition}
\newtheorem{Notation}[satz]{Notation}

\theoremstyle{remark}
\newtheorem{remark}[satz]{Remark}

\newcommand{\Z}{\mathbb{Z}}
\newcommand{\IZ}{\mathbb{Z}}
\newcommand{\C}{\mathcal{C}}

\newcommand{\zX}{X^{(*)}}
\newcommand{\zA}{A^{(*)}}

\newcommand{\variable}{\underline{\;\;}}

\providecommand{\Mf}[1]{\langle#1\rangle}

\DeclareMathOperator{\hocolim}{hocolim}
\DeclareMathOperator{\colim}{colim}
\DeclareMathOperator{\im}{im}
\DeclareMathOperator{\Tor}{Tor}
\DeclareMathOperator{\id}{id}

\DeclareMathOperator{\llcm}{l-lcm}

\DeclareMathOperator{\hgta}{ht_1}
\DeclareMathOperator{\hgtb}{ht_2}

\DeclareMathOperator{\Cat}{Cat}
\DeclareMathOperator{\sSet}{sSet}
\DeclareMathOperator{\Fun}{Fun}
\newcommand{\wt}{\widetilde}

\providecommand{\abs}[1]{\lvert#1\rvert}

\newenvironment{pf}{\begin{proof}[Proof]}{\end{proof}}

\begin{document}

\title[Discrete Morse Theory and Reformulation of the $K(\pi,1)$-conjecture]{Discrete Morse Theory and a Reformulation of the $K(\pi,1)$-conjecture}
\author{Viktoriya Ozornova}
\address{Institut für Mathematik, FU Berlin, Arnimallee 7, 14195 Berlin}
\thanks{ \texttt{viktoriya.ozornova@fu-berlin.de}}

\maketitle
\vspace{-1cm}
\begin{abstract}
A recent theorem of Dobrinskaya \cite{Dobrinskaya} states that the $K(\pi,1)$-conjecture holds for an Artin group $G$ if and only if the canonical map $BM\to BG$ is a homotopy equivalence, where $M$ denotes the Artin monoid associated to $G$. The aim of this paper is to give an alternative proof by means of discrete Morse theory and abstract homotopy theory. Moreover, we exhibit a new model for the classifying space of an Artin monoid, in the spirit of \cite{CMW}, and a small chain complex for computing its monoid homology, similar to the one of \cite{Squier}. 
\end{abstract}

\section{Introduction}
\label{Introduction}

Artin groups are a natural generalization of braid groups, intensively studied and well-understood in quite a few cases. Artin groups are closely related to Coxeter groups. The best understanding has been obtained in the case of Artin groups of finite type, i.e., those which correspond to finite Coxeter groups. Already in the early 70's, Brieskorn and Saito \cite{BrieskornSaito} as well as Deligne \cite{Deligne} solved the word and conjugacy problems for these groups and provided finite models for classifying spaces, which arise from the standard representation of the corresponding Coxeter groups. 

For Artin groups that are not of finite type, many questions are still open. It is in general not known whether they are torsion-free, have solvable word problem or a finite $K(\pi,1)$-model; for detailed accounts on these problems, we refer for example to \cite{CharneySurvey} and \cite{GodelleParisNeu}.

 One of the open conjectures about Artin groups is the so-called $K(\pi,1)$-conjecture. It says that a certain finite-dimensional space arising from the standard representation of the corresponding Coxeter group, similarly as in the finite-type case, is a model for the classifying space of the Artin group. This conjecture is known to be true for Artin groups of large type \cite{Hendriks}, for Artin groups of FC-type \cite{CD95} and for some groups of affine type \cite{Okonek}. There are also several reformulations of the problem, cf. e.g. \cite{CD95}, \cite{CD}, \cite{ParisKpi1}, \cite{Salvetti}. 

We are dealing with a rather new reformulation of this conjecture due to Dobrinskaya \cite{Dobrinskaya}. Each Artin group has an associated Artin monoid, introduced by Brieskorn and Saito \cite{BrieskornSaito}. They were also the first to use it to obtain information about the Artin group, transferring the solution of the word and conjugation problem from the monoid to the group for Artin monoids of finite type. In the finite type case, this is possible since then - and only then - Artin monoids satisfy Ore's condition. In the general case, it was even unclear for a long time whether Artin monoids inject into the corresponding Artin group. This was shown in 2002 by Paris \cite{ParisArtin}. For a monoid $M$ satisfying Ore's condition, it is always true that the map into the associated group $i\colon M\to G(M)$ induces a homotopy equivalence on classifying spaces. It is natural to ask whether this also holds for general Artin monoids. According to a result of Dobrinskaya, this question is equivalent to the $K(\pi,1)$-conjecture:
 \begin{Hauptsatz}[\cite{Dobrinskaya}] 
 The inclusion $BM \to BG$ is a homotopy equivalence if and only if the $K(\pi,1)$-conjecture holds for the Artin group $G$.
\end{Hauptsatz}

The first goal of this article is to give a new proof of this result. We use homotopy theory, in particular some results about homotopy colimits, instead of configuration spaces for the proof. The combinatorics are reminiscent of these by Dobrinskaya, yet they are arranged via the method of discrete Morse theory, which makes them more transparent. 

Discrete Morse theory is a tool to reduce a cellular object (e.g., simplicial complex, CW-complex, but also a based chain complex) to a homotopy equivalent, smaller one. It was proposed by Forman (cf. \cite{Forman}, see also \cite{FormanGuide}, \cite{Kozlov}, \cite{Chari} and similar ideas in \cite{Brown}) and developed in the last years to a useful tool in many areas. The idea is to give a coherent, combinatorial pattern for performing successive elementary collapses, which then provides a homotopy equivalence. We are going to use a version due to Batzies \cite{Batzies}, which is appropriate for infinite CW-complexes. 

Moreover, we use discrete Morse theory to exhibit a further model for a classifying space of an Artin monoid, which is a subcomplex of the bar complex model:
\begin{Hauptsatz}
Let $M$ be an Artin monoid, $\mathcal{E}$ a generating set closed with respect to left least common multiple and left complement. Then the subcomplex of $BM$ with cells given by 
\begin{eqnarray*}
 \mathcal{E}_*=\bigcup_{n\geq 0}\{[x_n|\ldots|x_1]\in BM\left| \mbox{ For all } 1\leq k \leq n, x_k\ldots x_1 \in \mathcal{E}\right\}.
\end{eqnarray*}
is homotopy equivalent to $BM$. In particular, there is a $\mathbb{Z}$-module complex computing the homology of $M$, with basis $\mathcal{E}_*$ as defined above and differentials given by restriction of the bar differential. 
\end{Hauptsatz}
 This is an analogue of the theorem by Charney, Meier and Whittlesey \cite{CMW} for Garside groups. For the details and the proof, see Section \ref{Matching on the bar complex}. Furthermore, we provide by means of discrete Morse theory a small chain complex computing the homology of an Artin monoid, having the same size as the one of Squier \cite{Squier}.

\subsection*{Overview}
\label{Overview}
The structure of this article is as follows. In Section \ref{Artin and Coxeter Groups}, we recall briefly the main properties of Artin groups which we will need later. In Section \ref{Discrete Morse Theory for Graded CW-complexes}, we describe a variant of discrete Morse theory for infinite CW-complexes due to Batzies \cite{Batzies}. In Section \ref{Geometric realization}, we make some preparations for applying discrete Morse theory to our specific situation. In Section \ref{A reformulation of the Kpi1-conjecture}, we present a new proof of Dobrinskaya's theorem (cf. \cite{Dobrinskaya}) using homotopy theory and discrete Morse theory. In Section \ref{Discrete Morse Theory for Chain Complexes}, we recollect another variant of discrete Morse theory, suitable for chain complexes. We apply it in Section \ref{Squier Complex for Artin Monoids} to obtain a chain complex which computes the homology of an Artin monoid and has the same size as Squier's complex \cite{Squier} and Salvetti's complex \cite{Salvetti}. In Section \ref{Matching on the bar complex}, we apply discrete Morse theory to obtain a small model for the classifying space of an Artin monoid, similar to the one by Charney, Meier and Whittlesey for Garside groups \cite{CMW}.

\subsection*{Acknowledgments}
\label{Acknowledgements}
The material of this note is a part of the author's PhD thesis at the University of Bonn. I would like to thank my advisor C.-F.~Bödigheimer for his supervision and support  during my PhD time. I also want to thank A.~Heß, who introduced me to discrete Morse theory, and who suggested an earlier version of Section \ref{Matching on the bar complex}, and Lennart Meier and C.-F.~Bödigheimer for proofreading the earlier versions of this paper. Furthermore, I would like to thank everybody who supported me during my PhD time or after it. In particular, I would like to thank the GRK 1150 Homotopy and Cohomology, the International Max Planck Research School on Moduli Spaces at MPIM Bonn and the SFB 647 Space-Time-Matter for their financial support.

\section{Artin and Coxeter Groups}
\label{Artin and Coxeter Groups}

We give a very brief review of Artin and Coxeter groups, mainly to fix the notation. There are many detailed accounts on these topics in the literature; see e.g. \cite{Bourbaki}, \cite{Humphreys}, \cite{ParisKpi1}, \cite{CharneySurvey}. 

\begin{Definition}\label{ExampleDefinitionArtinGroups}
An \textbf{Artin group} is a group given by a group presentation of the form
\begin{eqnarray*}
 G(S)=\Mf{S| \underbrace{sts\ldots}_{m_{s,t}} =\underbrace{tst\ldots }_{m_{s,t}}\mbox{ for all } s\neq t \in S}.
\end{eqnarray*}
where $m_{s,t}$ are natural numbers $\geq 2$ or infinity, with $m_{s,t}=m_{t,s}$ for all $s\neq t \in S$. Here, $m_{s,t}=\infty$ means that the pair $s,t$ does not satisfy any relation. We can associate to each Artin group a \textbf{Coxeter group} $W(S)$ by adding relations $s^2=1$ for all $s\in S$. (It is then consistent to set the numbers $m_{s,s}$ to be $1$.) The matrix $M_S=(m_{s,t})_{s,t\in S}$ will be called the \textbf{Coxeter matrix} defining $G(S)$ or $W(S)$, and the pair $(S, M_S)$ will be also called the \textbf{Coxeter system}. 

For each Coxeter system, we can define the corresponding \textbf{(positive) Artin monoid} by the monoid presentation
\begin{eqnarray*}
  M(S)= \Mf{S| \underbrace{sts\ldots}_{m_{s,t}} =\underbrace{tst\ldots }_{m_{s,t}}\mbox{ for all } s\neq t \in S}.
\end{eqnarray*}
For later use, we will denote the alternating word $sts\ldots$ with ${m}$ factors by $\Mf{s,t}^{m}$. Note that there are a monoid homomorphism $\pi\colon M(S)\to W(S)$ and a group homomorphism $\pi\colon G(S)\to W(S)$ mapping each generator to its image in the quotient group. 
We will call $M(S)$ as well as $G(S)$ or sometimes, by abuse of notation, even $S$ \textbf{of finite type} if the associated Coxeter group $W(S)$ is finite.
\end{Definition}

\begin{remark} 
  Artin monoids of finite type determine the behavior of the corresponding Artin group almost completely, as first shown by E.~Brieskorn and K.~Saito \cite{BrieskornSaito}. In the same article, they investigate more generally the structure of all Artin monoids. Amongst other things, they show that any Artin monoid is left and right cancellative.
  
  Observe moreover that, since the defining relations of an Artin monoid are length-preserving, there are no non-trivial invertible elements in an Artin monoid.
\end{remark}

We repeat briefly the basic notions for divisibility in monoids. They resemble divisibility in natural numbers, yet, one has to take into account that general monoids are non-commutative, so one has to distinguish between left and right divisibility. 

\begin{Definition}
Let $M$ be a monoid and let $x,y$ be elements in $M$. We say ``$x$ is a \textbf{left divisor} of $y$'' or, equivalently, ``$y$ is a \textbf{right multiple} of $x$'', and write $x\preceq y$ if there is an element $z\in M$ such that $y=xz$. Symmetrically, we define right divisors and left multiples; all notions for divisibility introduced later will have a symmetric analogue.

In a cancellative monoid without non-trivial invertible elements, the left and right divisibility relations are partial orders.
\end{Definition}

\begin{remark}
 A \textbf{left least common multiple} $c$ of two elements $a,b$ of a monoid $M$ is a left common multiple of these elements with the following property: Whenever $d$ is a left common multiple of $a$ and $b$, we have $d\succeq c$. This should not be confused with the notion of left minimal common multiple of $a$ and $b$, meaning a left common multiple of $a$ and $b$ which is not right-divisible by any other left common multiple of $a$ and $b$. 
 
 In general, least common multiples may or may not exist. 
  In a cancellative monoid without non-trivial invertible elements, left least common multiples are unique whenever they exist. 
  
  Inductively, one can also define the left least common multiple of any finite set of elements. They will be denoted by $\llcm$.
\end{remark}

For later use, we reassemble some results by Brieskorn and Saito \cite{BrieskornSaito}. For $I\subset S$, let $W(I)$ be the Coxeter group given by the restriction of the Coxeter matrix to $I$. (It is well-known that $W(I)$ is the subgroup of $W(S)$ generated by $I$). 
Recall that in \cite{BrieskornSaito}, it is shown that the Coxeter generating system $J$ of an Artin monoid $M(J)$ has a left least common multiple if and only if $W(J)$ is finite. We will denote this common multiple by $\Delta_{J}$ if it exists. Now let again $M(S)$ be any Artin monoid. We will consider the set of all such common multiples assigned to subsets of $S$:
\begin{eqnarray*}
\mathcal{D}=\{ \Delta_{J}\,|\, \varnothing \neq J\subset S, W(J)\mbox{ is finite} \}. 
\end{eqnarray*}
This subset of $M$ is clearly a generating set, since for each $s\in S$, the group $W(\{s\})\cong \mathbb{Z}/2$ is finite and thus $\Delta_{\{s\}}=s\in \mathcal{D}$. 

Brieskorn and Saito exhibit normal forms for elements of $M(S)$ with respect to $\mathcal{D}$. Although these normal forms are not geodesic, they are a helpful tool in studying Artin monoids. 

For $x\in M$, let $I(x):=\{a\in S|\; \exists \,y\in M\colon \, x=ya\}$ be the set of letters in $S$ with which a word for $x$ may start (on the right). Note that $I(\Delta_J)=J$ for any subset $J\subset S$ admitting a left least common multiple (cf. e.g. \cite{Squier}, Part II, §5). The Brieskorn-Saito normal form is given by the following theorem.

\begin{Theorem}[\cite{BrieskornSaito}, §6]\label{BSNF}
 For any $w\in M$, there are unique non-empty subsets 
\begin{eqnarray*}
 I_k, \ldots, I_1 \subset S
\end{eqnarray*}
 such that 
\begin{eqnarray*}
 w=\Delta_{I_k}\Delta_{I_{k-1}}\ldots \Delta_{I_2}\Delta_{I_1}
\end{eqnarray*}
and $I(\Delta_{I_k}\ldots \Delta_{I_j})=I_j$ for $1\leq j\leq k$.
\end{Theorem}

\section{Discrete Morse Theory for Graded CW-complexes}
\label{Discrete Morse Theory for Graded CW-complexes}

We start by describing a version of discrete Morse theory due to E.~Batzies \cite{Batzies}. We follow his exposition very closely. This version can be applied to infinite CW-complexes; furthermore, Batzies formulates the theory in the language of acyclic matchings (instead of discrete Morse functions), which seems to be convenient for our approach. The approach by E.~Batzies generalizes the original one due to R.~Forman \cite{Forman}, and of M.~Chari \cite{Chari} who was first to work with acyclic matchings instead of Morse functions. K.~Brown \cite{Brown} uses a variant of discrete Morse theory for simplicial sets similar to the one we will use. We stick to the version by Batzies since this gives a more detailed description of the Morse complex.  We start with the basic definitions. Instead of the partially ordered set of simplices in a simplicial complex, we will consider the poset of cells in the CW complex, defined as follows.

\begin{Definition}
Let $X$ be a CW-complex, and let $X^{(*)}$ be the set of its open cells. For two cells $\sigma, \sigma'\in \zX$, we write $\sigma \leq \sigma'$ iff the closed cell $\overline{\sigma}$ is a subset of a closed cell $\overline{\sigma'}$, and call $\sigma$ a \textbf{face} of $\sigma'$. We say that a cell $\sigma$ is a \textbf{facet} of a cell $\sigma'$ if $\sigma\neq \sigma'$, $\sigma\leq \sigma'$ and for any $\tau \in \zX$ with $\sigma\leq \tau\leq \sigma'$ we have either $\sigma=\tau$ or $\sigma'=\tau$. 

If $(P, \preceq)$ is any poset, a $P$-\textbf{grading} on $X$ is a poset map $f\colon \zX \to P$. Given a $P$-grading $f$ and $p\in P$, we write $X_{\preceq p}$ for the sub-CW-complex of $X$ consisting of all cells $\sigma$ with $f(\sigma)\preceq p$.
\end{Definition}

To perform elementary collapses given by the Morse equivalence (to be defined below), we have to ensure some regularity for the cells in question.

\begin{Definition}[\cite{Forman}]
Let $X$ be a CW complex, let $\sigma$ be an $n$-dimensional cell of $X$ and let $\tau$ be an $(n+1)$-dimensional cell with characteristic map $f_{\tau}\colon D^{n+1}\to X$. Assume $\sigma\leq \tau$. We call $\sigma$ a \textbf{regular face} of $\tau$ if $f_{\tau}$ restricted to $f_{\tau}^{-1}(\sigma)$ is a homeomorphism onto $\sigma$ and, in addition, $\overline{f_{\tau}^{-1}(\sigma)}$ is a closed metric $n$-ball in $\partial D^{n+1}$.
\end{Definition}

\begin{remark}
Here, for the definition of regular faces, we use the terminology by Forman \cite{Forman}. This is compatible with all the proofs in \cite{Batzies}. 
\end{remark}

The combinatorial data we will use is the one of an acyclic matching, similar to the notion of the gradient vector field of a usual Morse function. 

\begin{Definition} \label{BatziesAcyclicMatching}
 Let $X$ be a CW-complex. The \textbf{cell graph} $G_X$ of $X$ is a directed graph with $\zX$ as the set of vertices. There is an edge from $\sigma$ to $\tau$ (denoted by $\sigma \to \tau$) if and only if $\tau$ is a facet of $\sigma$. In other words, the set of edges is given by
\begin{eqnarray*}
 E_X:=\{\sigma\to \tau| \tau \mbox{ is a facet of }\sigma\}.
\end{eqnarray*}
A \textbf{matching} on $X$ is a subset $A\subset E_X$ such that the following conditions hold:
\renewcommand{\labelenumi}{(M\arabic{enumi})}
\begin{enumerate}
 \item If $(\sigma \to \tau)\in A$, then $\tau$ is a regular face of $\sigma$. 
\item Each cell of $X$ occurs in at most one edge of $A$.
\end{enumerate}
We associate to a matching $A$ a new graph $G_X^A$ by inverting all arrows in $A$ and keeping all other arrows unchanged. More precisely, $G_X^A$ is a directed graph with the same vertices as $G_X$ and with edge set 
\begin{eqnarray*}
 E_X^A:=\left (E_X\setminus A\right )\cup\{\sigma \to \tau|(\tau\to \sigma)\in A\}.
\end{eqnarray*}

A matching $A$ is called \textbf{acyclic} if in addition we  have
\begin{enumerate}
\setcounter{enumi}{2}
\item The graph $G_X^A$ contains no cycle.
\end{enumerate}
A cell of $X$ is called $A$-\textbf{essential} if it does not occur in $A$. We denote by $\zX_{ess}$ the set of essential cells of $X$.

If $(\sigma \to \tau)$ is an element of $A$, we call the cell $\tau$ \textbf{redundant} and the cell $\sigma$ its \textbf{collapsible} partner.
\end{Definition}

Such a matching defines now a new poset as follows. 

\begin{Definition}
Let $X$ be a CW-complex and $A$ an acyclic matching on it. 
 We set $\zA=A\sqcup \zX_{ess}$ as sets, so an element in $\zA$ is either an essential cell or an edge belonging to the matching. Now we define a partial order on $\zA$ as follows: Let $\widetilde{G}_X^{A}$ be a graph with vertices $\zX$ and edge set 
\begin{eqnarray*}
 E_X^A:=E_X\cup\{\sigma \to \tau|(\tau\to \sigma)\in A\},
\end{eqnarray*}
i.e., we add to $G_X$ all reversed edges of $A$. For $a,b\in \zA$, we set $a\preceq_A b$ if there is a path in $\widetilde{G}_X^{A}$ from $b$ to $a$. If $b$ is an element of the form $\sigma \to \tau$, this means that the path may start either from $\sigma$ or from $\tau$; similarly, if $a$ is of the form $\sigma \to \tau$, the path may end either at $\sigma$ or at $\tau$. 

This defines a partial order on $\zA$. We call the poset $(\zA, \preceq_A)$ the \textbf{matching poset} of $A$. 

The map given by
\begin{eqnarray*}
 \zX & \to  &\zA\\
 \sigma & \mapsto & \begin{cases}
                \sigma, \mbox{ if } \sigma\in \zX_{ess},\\
		(\tau\to \tau'), \mbox{ if } \sigma \in \{\tau, \tau'\} \mbox{ and } (\tau\to \tau')\in A
               \end{cases}
\end{eqnarray*}
can be seen to be order-preserving. We call it the \textbf{universal $A$-grading} on $X$.
\end{Definition}

We will need some finiteness conditions to handle our CW complexes, which are often not finite dimensional. 

\begin{Definition}
 Let $(P,\preceq)$ be a poset and $f\colon \zX\to P$ a grading on a CW complex $X$. We call the grading $f$ \textbf{compact} if $X_{\preceq p}$ is compact for all $p\in P$. 
\end{Definition}

Last, we need the definition of the Morse complex of a matching. It is quite technical. We will still cite it here since we will need it quite explicitly. The intuition behind it is to glue a cell to the new cell complex for each essential cell of the old one. Since an essential cell may have had redundant or collapsible cells in its boundary, we have to change the gluing maps appropriately. 

\begin{Definition} \label{MorseComplexofMatching}
 Let $X$ be a CW complex and $A$ an acyclic matching on it such that the universal $A$-grading is compact. For all $a\in \zA$, we define first inductively $(X_A)_{\preceq a}$, and also a map $H(A)_{\preceq a}\colon X_{\preceq a} \to (X_A)_{\preceq a}$. In the end, these will be the pieces of the Morse complex and of the homotopy equivalence from $X$ to the Morse complex. 

First, if $a\in\zA$ is minimal, we know that $a\in \zX_{ess}$ and $X_{\preceq a}=a$. We define $(X_A)_{\preceq a}$ to be equal to $a$ and the map $H(A)_{\preceq a}$ to be just the identity. 

Now take any $a\in\zA$ and suppose the associated piece of the Morse complex $(X_A)_{\preceq b}$ and the Morse equivalence $H(A)_{\preceq b}$ are already constructed for all $b\prec a$ in a way such that $b\preceq b'\prec a$ induces an inclusion of the associated pieces of the Morse complex and the restrictions of the future Morse equivalence are compatible with these inclusions. Then define first
\begin{eqnarray*}
 (X_A)_{\prec a}:=\bigcup_{b\prec a} (X_A)_{\preceq b}
\end{eqnarray*}
to be the colimit over the poset $\{b\prec a\}$ of already known pieces and let the map $H(A)_{\prec a}$ from $X_{\prec a}$ be induced by the already known pieces. Now we have to distinguish whether $a$ is an element of $A$ or of $\zX_{ess}$. If $a=(\tau \to \sigma)\in A$, then we define $(X_A)_{\preceq a}=(X_A)_{\prec a}$ and let the map be defined by 
\begin{eqnarray*}
 H(A)_{\preceq a}=H(A)_{\prec a}\circ \widetilde{h}_{\tau\to\sigma}
\end{eqnarray*}
where the map $\widetilde{h}_{\tau\to\sigma}$ deforms $X_{\preceq a}$ into $X_{\prec a}$ by deforming $\overline{\tau}$ into the union of its faces different from $\sigma$. This is possible since $\sigma$ is a regular face of $\tau$; for more details, we refer again to \cite{Batzies}. 

Now we consider the other case $a=\sigma \in \zX_{ess}$, where $\sigma$ is a cell of dimension $i$ with characteristic map $f_{\sigma}\colon D^i\to X_{\preceq a}$. We define 
\begin{eqnarray*}
 (X_A)_{\preceq a}=D^i\cup_{H(A)_{\prec a}\circ f_{\partial \sigma}} (X_A)_{\prec a}
\end{eqnarray*}
so we glue a new cell to $(X_A)_{\prec a}$ via $H(A)_{\prec a}\circ f_{\partial \sigma}$. The new piece of map is now induced by the identity on the new cell: Define
\begin{eqnarray*}
 H(A)_{\preceq a}=\id_{D^i} \cup_{f_{\partial \sigma}} H(A)_{\prec a}. 
\end{eqnarray*}

Last, define the \textbf{Morse complex} $X_A$ to be the colimit of all pieces and the \textbf{Morse equivalence} $H(A)\colon X\to X_A$ to be the induced map on it. 
\end{Definition}

We will need the following theorem which is a version of the main theorem of discrete Morse theory in Batzies' flavor. 
\begin{Theorem}(\cite{Batzies}) \label{gradDMTMainThm}
 Let $X$ be a CW complex and $A$ an acyclic matching on it such that the universal $A$-grading is compact. Then the $i$-cells of the Morse complex $X_A$ are in one-to-one correspondence with the essential cells of $A$ of dimension $i$. Furthermore, the Morse equivalence $H(A)\colon X \to X_A$ is a homotopy equivalence. 
\end{Theorem}

Last, we will need a criterion to check whether the universal $A$-grading is compact. We will use the following lemma.

\begin{Lemma}(\cite{Batzies})\label{LemmaCompact}
Let $X$ be a CW complex and $A$ an acyclic matching on it. Furthermore, let $P$ be a poset and let $f\colon \zX\to P$ be a compact grading on $X$ such that $f(\tau)=f(\sigma)$ holds for all $(\tau\to \sigma)\in A$. Then the universal $A$-grading is also compact.
\end{Lemma}

We derive a corollary of Theorem \ref{gradDMTMainThm}. 
\begin{Cor} \label{UnterkomplexDMT}
Let $X$ be a CW complex and $A$ an acyclic matching on it such that the universal $A$-grading is compact. Assume furthermore that the essential cells of $A$ form a subcomplex $X_{ess}$ of $X$, i.e., if $\sigma \in \zX_{ess}$ and $\tau\leq \sigma$, then $ \tau \in \zX_{ess}$. Then the inclusion $i\colon X_{ess}\to X$ is a homotopy equivalence. 
\end{Cor}

\begin{pf}
We will show that the composition $H(A)\circ i\colon X_{ess} \to X_A$ is a homotopy equivalence; this will imply the claim. More precisely, we will first show inductively that $(X_{ess})_{\preceq a}=(X_A)_{\preceq a}$ for all $a\in\zA$ and the map $H(A)\circ i$ is the identity. For $a\in \zA$ minimal, the statement is clear. Assume we have proven the statement for all $b\prec a$ and we would like to show it for $a$. If $a$ is of the form $(\tau\to \sigma)$ in $A$, then $(X_{ess})_{\preceq a}=(X_{ess})_{\prec a} \subset X_{\prec a}$. Note that $\widetilde{h}_{\tau\to\sigma}$ is identity on $(X_{ess})_{\preceq a}$, so that we are done in this case.

Now assume that $a=\sigma$ is an essential cell of dimension $n$. Let $f_{\sigma}\colon D^n\to X_{\preceq a}$ be the characteristic map of this cell. Note that by assumption the attaching map $f_{\partial\sigma}$ has its image in $(X_{ess})_{\prec a}$ and $(X_{ess})_{\preceq a}=(X_{ess})_{\prec a}\cup_{f_{\partial\sigma}} D^n$. It is also $(X_A)_{\preceq a}$ by the induction hypothesis and by the definition of $X_A$. Moreover, the composition of the inclusion with $H(A)_{\preceq a}$ is again the identity. This completes the induction step. (Observe that the compactness of the grading enables the induction arguments.) Taking the union of all $(X_A)_{\preceq a}$, we see that $X_{ess}=X_A$ and $H(A)\circ i$ is the identity.

Altogether, we have shown that $X_{ess}\to X$ is a homotopy equivalence.
\end{pf}

\section{Geometric Realization}
\label{Geometric realization}

The aim of this section is to make some observations about properties of geometric realization which will be used later on. We start with the following two well-known properties:

\begin{Prop}[\cite{FritschPiccinini}, Section 4.3]\label{RealisierungsLemma}
 \begin{enumerate}
  \item The geometric realization of a simplicial set $X$ is a quotient space of the subspace $\coprod X_n^{\#}\times \Delta^n$ of $\coprod X_n\times \Delta^n$, where $X_n^{\#}$ denotes the set of non-degenerate $n$-dimensional simplices of $X$. 
 \item For a simplicial set $X$, each point of the geometric realization $\abs{X}$ has a unique presentation as a pair $(x, \underline{t})$, where $x$ is non-degenerate and $\underline{t}\in\Delta^{\dim x}$ is an inner point. 
 \end{enumerate}

\end{Prop}

We prove the following easy consequence:
\begin{Lemma}\label{RealisierungEigenschaftA}
 Let $X$ be a simplicial set with the following property: All faces of a non-degenerate simplex are again non-degenerate. Then there is a homeomorphism 
\begin{eqnarray*}
r(X):=  \left(\coprod X_n^{\#}\times \Delta^n/\sim\right ) \to \abs{X}
\end{eqnarray*}
where $\sim $ is generated by $(d_ix, (t_0,\ldots, t_{n-1}))\sim (x, (t_0, \ldots, t_{i-1}, 0, t_i, \ldots, t_{n-1}))$. Moreover, the projection 

\begin{eqnarray*}
\coprod X_n^{\#}\times \Delta^n \to \coprod X_n^{\#}\times \Delta^n/\sim 
\end{eqnarray*}
defines a CW structure on $r(X)$. Furthermore, each element of $r(X)$ has a unique representative of the form $(x, \underline{t})$, where $\underline{t}\in\Delta^{\dim x}$ is an inner point. 
\end{Lemma}

\begin{pf}
 Since both $r(X)$ and $\abs{X}$ are quotient spaces of $\coprod X_n^{\#}\times \Delta^n$ (using Proposition \ref{RealisierungsLemma}), it is enough to construct mutually inverse bijections $r(X)\to \abs{X}$ and $\abs{X}\to r(X)$ which are compatible with the quotient maps. Then, by the definition of the quotient topology, both maps are continuous and thus homeomorphisms. 

The map $f\colon r(X)\to \abs{X}$ is given by simply regarding an equivalence class $[x, \underline{t}]$ in $r(X)$ as an equivalence class in $\abs{X}$. This is clearly well-defined and compatible with the quotient maps. 

For the other direction, we take any $[y,\underline{s}]\in \abs{X}$ and consider its unique representative $[x, \underline{t}]$ as in Proposition \ref{RealisierungsLemma}. Since $x \in X_m^{\#}$ for some $m$, it also defines a point $g([y,\underline{s}])$ in $r(X)$. This gives us again a well-defined map, which is obviously compatible with the quotient maps.

It is also immediate that $fg=\id$. For the other direction, let $[x, \underline{t}]\in r(X)$ and assume $\underline{t}$ is not an inner point of $\Delta^{\dim x}$. Then there is an inner point $\underline{u}\in \Delta^{m}$ and a sequence of natural numbers $i_1, \ldots, i_k$ such that $\underline{t}=\delta_{i_1}\ldots \delta_{i_k}(\underline{u})$. Then 
\begin{eqnarray*}
 [x, \underline{t}]=[d_{i_k}\ldots d_{i_1}(x), \underline{u}] \in r(X)
\end{eqnarray*}
where $d_{i_k}\ldots d_{i_1}(x)$ is again a non-degenerate simplex by assumption. This shows that also $gf=\id$. Altogether, this proves the first claim. 

The second claim is completely analogous to the statement that $\abs{X}$ is a CW complex. 

The last claim follows immediately from the second part of the Proposition \ref{RealisierungsLemma}. 
\end{pf}

\begin{remark} \label{RealisierungR}\hspace*{\fill} \\\vspace*{-0.5cm} 
\begin{enumerate}
 \item Simplicial sets as in Lemma \ref{RealisierungEigenschaftA} are said to have \textbf{Property A}, e.g., in \cite{MayNotes}. 
 \item We will from now on identify $r(X)$ and $\abs{X}$ under the conditions of the last lemma since these spaces are then homeomorphic and have the ``same'' CW structure.
 \end{enumerate}
\end{remark}

In discrete Morse theory, we have to check whether a smaller cell is a regular face of a larger one. We provide for this purpose a regularity criterion for realizations of simplicial sets. 

\begin{Lemma}\label{RegularityCriterion}
Let $Y$ be a simplicial set fulfilling Property A and let $s$ be a non-degenerate $n$-simplex in $Y$. Consider $t=d_i(s)$ for some $0\leq i\leq n$. Let $\sigma$ and $\tau$ be cells of $r(Y)$ (as defined above) corresponding to $s$ and $t$, respectively. If $d_j(s)\neq t$ for all $0\leq j\neq i\leq n$, then $\tau$ is a regular face of $\sigma$.
\end{Lemma}

\begin{pf}
Fix a homeomorphism $\psi\colon D^n\to \Delta^n$ such that $S^{n-1}_{\geq 0}$ is mapped homeomorphically to the $i$-th side of $\Delta^n$. The map 
\begin{eqnarray*}
f_{\sigma}\colon\{s\}\times D^n\xrightarrow{\psi}\{s\}\times\Delta^n\hookrightarrow \coprod_j Y_j^{\#}\times \Delta^j \twoheadrightarrow r(Y)
\end{eqnarray*}
is the characteristic map of $\sigma$. Any point in the (open) cell $\tau$ is of the form 
\begin{eqnarray*}
x=[t, (t_0, \ldots, t_{n-1})]
\end{eqnarray*}
with $t_i>0$ and $(t_0, \ldots, t_{n-1})\in \Delta^n$. Note that this is also the unique representative with an inner point in the second coordinate, as described in Lemma \ref{RealisierungEigenschaftA}. This point is by definition identified with the point represented by 
\begin{eqnarray*}
 (s, (t_0, \ldots, t_{i-1}, 0, t_i, \ldots, t_{n-1})) 
\end{eqnarray*}
of $\{s\}\times \Delta^n$. Assume there is another point with representative of the form $(s, \underline{u})$ which is identified with $x$. Then $\underline{u}$ cannot be an inner point by the uniqueness statement of Lemma \ref{RealisierungEigenschaftA}. So we can write $\underline{u}=\delta_{i_1}\ldots \delta_{i_k}(\underline{v})$, where $k\geq 1$ and $\underline{v}$ is an inner point of an appropriate simplex. Thus, $x$ has also a representative of the form $(d_{i_k}\ldots d_{i_1}(s), \underline{v})$. Using again the uniqueness, we see that $\underline{v}=(t_0,\ldots, t_{n-1})$ and $d_{i_k}\ldots d_{i_1}(s)=t$. This implies that $k=1$. By hypothesis of the lemma, $d_{i_1}(s)=t$ implies $i_1=i$. This implies that $f_{\sigma}$ is injective when restricted to $f_{\sigma}^{-1}(\tau)$, where the last one is the interior of $\Delta^{n-1}$ considered as $i$-th boundary of $\Delta^n$. Thus, the second condition for regularity is already fulfilled. Furthermore, the map $f_{\sigma}\colon D^n \to \overline{\sigma}$ is an identification. It is a simple observation that for an identification  map $q\colon Z\to Z'$ and $B\subset Z'$ open or closed subset, the restriction $q\colon q^{-1}B\to B$ is an identification again (cf. e.g. the textbook by T.~tom Dieck \cite{tomDieck}). Thus the restriction of $f_{\sigma}$ to $f_{\sigma}^{-1}(\overline{\tau})$ is an identification since $\overline{\tau}\subset \overline{\sigma}$ is closed. We can now apply same argument again since $\tau \subset \overline{\tau}$ is open in $\overline{\tau}$. This completes the proof that $f_{\sigma}$ restricted to $f_{\sigma}^{-1}(\tau)$ is a homeomorphism. 
\end{pf}

We will  later need the following easy lemma to apply our regularity criterion. For a small category $\mathcal{C}$, we denote by $N\mathcal{C}$ its nerve. 
\begin{Lemma}[\cite{MayNotes}, Lemma 11]\label{MayCriterion}
 Let $\mathcal{C}$ be a small category. Then $N\mathcal{C}$ has Property A if and only if the following holds: Whenever $f\colon A\to B$ and $g\colon B\to A$ are morphisms in $\mathcal{C}$ such that $g\circ f=\id_A$, then we already have $A=B$ and $f=g=\id$. 
\end{Lemma}

\section{A Reformulation of the $K(\pi, 1)$-conjecture}
\label{A reformulation of the Kpi1-conjecture}

The aim of this section is to reprove a theorem by N.~Dobrinskaya \cite{Dobrinskaya} claiming that the $K(\pi,1)$-conjecture for an Artin group $G(S)$ is equivalent to the statement that the inclusion $BM(S)\to BG(S)$ is a homotopy equivalence. Her proof is long and uses the machinery of configuration spaces. It seems that our proof is, though in a sense less geometrical, yet more transparent and less involved. The combinatorics entering is rather similar, but it seems to be more systematic to arrange them via discrete Morse theory.  

We will use a reformulation of the $K(\pi, 1)$-conjecture by R.~Charney and M.~Davis \cite{CD95}, reformulated in the language of Grothendieck constructions (see Conjecture \ref{Kp1ConjCD}). First, we introduce the necessary vocabulary and basic facts. 

Let $M:=M(S)$ be an Artin monoid with Artin-Coxeter generating set $S$ and let $G(S)$ be the corresponding (Artin) group and $W(S)$ the corresponding Coxeter group. 
 For $I\subset S$, let $W(I)$ be the Coxeter group given by the restriction of the Coxeter matrix to $I$. Let $S^f$ be the collection of all subsets $I$ of $S$ such that the Coxeter group $W(I)$ is finite.

\begin{Theorem}[\cite{CD95}]
The canonical map  $\colim_{T \in S^f} G(T) \to G(S)$ is an isomorphism. 
\end{Theorem}

We want to observe that the same holds for monoids.

\begin{Lemma}
  The canonical map $\colim_{T \in S^f} M(T)\to M(S)$ is an isomorphism.
\end{Lemma}

\begin{pf}
Let $N$ be any monoid, and assume we have compatible monoid homomorphisms $\varphi_{T}\colon M(T)\to N$ for $T\in S^{f}$. Since for each $a\in S$, we have $\{a\}\in S^{f}$, we can define a map $\varphi \colon S \to N$ via $\varphi(a):=\varphi_{\{a\}}(a)$. We now want to show that $\varphi$ defines a monoid homomorphism $M(S)\to N$. Recall that all relations in $M(S)$ are of the type 
\begin{eqnarray*}
\Mf{a,b}^{m_{a,b}}=\Mf{b,a}^{m_{a,b}}
\end{eqnarray*}
whenever $m_{a,b}$ is finite. But if $m_{a,b}$ is finite, the corresponding dihedral group associated to $M(\{a,b\})$ is finite, thus $\{a,b\}\in S^f$. Since $\varphi_{\{a,b\}}(a)=\varphi_{\{a\}}(a)=\varphi(a)$ and similar for $b$, we know that the elements $\varphi(a), \varphi(b)\in N$ satisfy
\begin{eqnarray*}
\Mf{\varphi(a),\varphi(b)}^{m_{a,b}}=\Mf{\varphi(b),\varphi(a)}^{m_{a,b}}
\end{eqnarray*}
since $\varphi_{\{a,b\}}$ is a monoid homomorphism. So $\varphi$ is a well-defined monoid homomorphism, and it is compatible with each $\varphi_T$ since they coincide on $T$. Moreover, since the values of $\varphi$ on $S$ are fixed by the family $\varphi_T$, the monoid homomorphism $\varphi$ is unique. So $M(S)$ has the universal property of the colimit, and this implies the claim.
\end{pf}

We will need the Grothendieck construction for a functor $F\colon \mathcal{C} \to \Cat$, where $\C$ is a small category and $\Cat$ is the category of small categories, as described by Thomason in \cite{Thomason}. The following definition is taken from \cite{Thomason}. 
To a functor $F\colon \mathcal{C} \to \Cat$ as above, we assign a category $\mathcal{C}\int F$, called the Grothendieck construction. Its objects are pairs $(C, x)$, where $C$ is an object in $\mathcal{C}$ and $x$ is an object in $F(C)$. A morphism from $(C_1,x_1)$ to $(C_2, x_2)$ is given by a map $c\colon C_1\to C_2$ in $\mathcal{C}$ and a map $\varphi\colon F(c)(x_1)\to x_2$ in the category $F(C_2)$. The composition with a further morphism $(c', \varphi')\colon (C_0, x_0)\to (C_1,x_1)$ is given by 
\begin{eqnarray*}
(c, \varphi)\circ (c', \varphi')=(cc', \varphi\circ F(c)(\varphi')).
\end{eqnarray*}

Note that the construction is functorial: a natural transformation $\alpha\colon F\Rightarrow F'$ induces a functor $\mathcal{C}\int\alpha\colon \mathcal{C}\int F\to \mathcal{C}\int F'$, given by $\left (\mathcal{C}\int\alpha\right )(C,x)=(C,\alpha(C)(x))$ on objects and by $\left (\mathcal{C}\int\alpha\right )(c,\varphi)=(c, \alpha(C_2)(\varphi))$ on morphisms. One checks that this defines a functor from the functor category $\Fun(\mathcal{C}, \Cat)$ into $\Cat$.

In \cite{Thomason}, Thomason identifies the nerve of the Grothendieck construction with a certain homotopy colimit in simplicial sets.  For the exact definition of a homotopy colimit, see e.g. \cite{BousfieldKan}. We will need mainly the following lemma and the theorem below:
\begin{Lemma}[\cite{BousfieldKan}, XII, 3.7 and 4.2]
Let $X,Y$ be two functors from a small category $\C$ to simplicial sets, and let $\psi\colon X\to Y$ be a natural transformation. Then there exists an induced map $\hocolim \psi\colon \hocolim X\to \hocolim Y$ making $\hocolim$ into a functor from the functor category $\Fun(\mathcal{C}, \sSet)$ into $\sSet$. Furthermore, if for all objects $C$ in $\mathcal{C}$, the map $\psi(C)\colon X(C)\to Y(C)$ is a weak homotopy equivalence, then the induced map $\hocolim \psi$ is also a weak homotopy equivalence.  
\end{Lemma}

We will use the following homotopy colimit theorem by Thomason:
\begin{Theorem}[\cite{Thomason}] \label{ThomasonTheorem}
Let $F\colon \mathcal{C}\to \Cat$ be a functor. Then there is a natural weak homotopy equivalence
\begin{eqnarray*}
\eta(F)\colon \hocolim NF \to N(\mathcal{C}\int F)
\end{eqnarray*}
of simplicial sets.
\end{Theorem}

Combining these two results, we obtain: 
\begin{Prop}
Let $F,G\colon \mathcal{C} \to \Cat$ be two functors starting from a small category $\mathcal{C}$, and let $\psi\colon F\Rightarrow G$ be a natural transformation between them such that 
\begin{eqnarray*}
 N(\psi(C))\colon NF(C)\to NG(C)
\end{eqnarray*}
 is a weak homotopy equivalence for each object $C$ of $\mathcal{C}$. Then the induced map of simplicial sets $N(\mathcal{C}\int \psi)\colon N(\mathcal{C}\int F)\to N(\mathcal{C}\int G)$ is a weak homotopy equivalence.
\end{Prop}

We will now apply this proposition to our situation. We consider the functors 
\begin{eqnarray*}
 M(-), G(-)\colon S^f\to \Cat
\end{eqnarray*}
 associating to $T\in S^f$ the corresponding Artin monoids and Artin groups, respectively. Here, the category $S^f$ means the category associated to the poset $S^f$ with usual inclusion as ordering, and monoids and groups are viewed as categories with one object. There is a natural transformation $i\colon M(-)\to G(-)$ given by the canonical map. By \cite{BrieskornSaito}, the Artin monoids of finite type satisfy the Ore condition and are cancellative. By \cite{CartanEilenberg}, Ch. X, §4, we know that for a cancellative monoid $M$ satisfying the Ore condition and its associated group $G$, the $\Tor$-term $\Tor^{\mathbb{Z}[M]}_n(\Z[G], \Z)$ vanishes for all $n>0$. Now we will use the following proposition of Fiedorowicz:

\begin{Prop}[\cite{Fiedorowicz}]\label{BMBGKrit}
 Let $M$ be a monoid and let $G$ be its associated group. Then, the following are equivalent:
\begin{enumerate}
 \item $\pi_k(BM)=0$ for all $k\geq 2$.
 \item The map $BM\to BG$ is a homotopy equivalence.
 \item $\Tor^{\mathbb{Z}[M]}_n(\Z, \Z[G])=0$ for all $n>0$.
\end{enumerate}
\end{Prop}

So we know that the inclusion $BM(T)\to BG(T)$ is a homotopy equivalence for $T\in S^f$. Altogether, we have proven:
\begin{Cor}\label{GrothendieckVergleich}
 The map  $N(S^f\int M(-))\to N(S^f\int G(-))$ induced by the inclusion $i$ is a weak homotopy equivalence. 
\end{Cor}

Now we are going to describe $S^f\int M(-)$ and $S^f\int G(-)$ more concretely. Since each $M(T)$ and each $G(T)$ has exactly one object, the set of objects of either Grothendieck construction is exactly $S^f$. There can be a map from $T$ to $T'$ only if $T\subset T'$. Each such map is given by a self-map of the only object of $M(T')$, so we have
\begin{eqnarray*}
S^f\int M(-)(T,T')=\begin{cases}
M(T'), \mbox{ if } T\subset T'\\
\varnothing, \mbox{ else.}
\end{cases}
\end{eqnarray*}
Note that the composition is given by the monoid multiplication. The category $S^f\int G(-)$ has a completely analogous description. We are now going to show:

\begin{Prop}\label{MainIngridientDobr}
The space $BM(S)$ is homotopy equivalent to $\abs{N(S^f\int G(-))}$. 
\end{Prop}

Before proving the Proposition, we will point out why this shows the desired equivalence. It follows from \cite{CD95}, Corollary 3.2.4, that the $K(\pi,1)$-conjecture for $G(S)$ is equivalent to the following:
\begin{Conj}[\cite{CD95}] \label{Kp1ConjCD}
 For any Artin group $G(S)$, the space $\abs{N(S^f\int G(-))}$ is homotopy equivalent to $BG(S)$. 
\end{Conj}

We can now use the results stated above to conclude the following theorem, first proven (by different means) by N.~Dobrinskaya. 

\begin{Theorem}(cf. \cite{Dobrinskaya}) \label{DobrinskayaThm}
 The inclusion $BM(S) \to BG(S)$ is a homotopy equivalence if and only if the space $\abs{N(S^f\int G(-))}$ is homotopy equivalent to $BG(S)$. 
\end{Theorem}

\begin{pf}
 If the inclusion $BM(S) \to BG(S)$ is a homotopy equivalence, Conjecture \ref{Kp1ConjCD} holds for $G(S)$ by Proposition \ref{MainIngridientDobr}. 

For the other implication, we use again Proposition \ref{BMBGKrit}: If the Conjecture \ref{Kp1ConjCD} holds for $G(S)$, again Proposition \ref{MainIngridientDobr} implies that all higher homotopy groups of $BM$ vanish, thus the claim. 
\end{pf}

We still need to show Proposition \ref{MainIngridientDobr}. First, we will reformulate the statement once more. We want to describe the space $\abs{N(S^f\int G(-))}$ differently so we can see that it is homotopy equivalent to $BM(S)$. Recall that by Corollary \ref{GrothendieckVergleich} and Theorem \ref{ThomasonTheorem}, we have identified the former space up to homotopy with $\abs{\hocolim_{S^f} NM(-)}$. We proceed by describing this space differently. Before doing so, we need an auxiliary result.

\begin{Lemma} \label{ColimitSimplicialSet}
 The simplicial set $\wt{K}:=\bigcup_{I\in S^f} NM(I)$ is the colimit of the functor $S^f \to \sSet$, given by $J\mapsto NM(J)$. 
\end{Lemma}

\begin{pf}
 Since there are compatible maps $NM(J)\to \wt{K}$, given just by inclusion, we obtain a map of simplicial sets
\begin{eqnarray*}
 \colim_{J\in S^f} NM(J) \to \wt{K}=\bigcup_{I\in S^f} NM(I).
\end{eqnarray*}
This map is obviously surjective. We will now show that it is also injective. Assume there are some $[x], [y] \in \colim_{J\in S^f} NM(J)$ coming from simplices $x\in NM(J_1)_k$, $y\in NM(J_2)_k$ and mapped to the same element in $\bigcup_{I\in S^f} NM(I)$. This implies that there is a simplex $z \in NM(J_1)_k\cap NM(J_2)_k=NM(J_1\cap J_2)_k$ mapping both to $x$ and to $y$ under corresponding inclusions. This implies exactly that $[x]=[z]=[y]\in \colim_{J\in S^f} NM(J)$, proving the injectivity.
\end{pf}

Now we are ready to show the following proposition.

\begin{Prop}\label{ModellKatHocolim}
 The functor $NM(-)\colon S^f \to \sSet$ is cofibrant as an object of the category $\Fun(S^f, \sSet)$, where the model structure on the latter is given by levelwise weak equivalences and levelwise fibrations. This implies in particular that $\wt{K}$ from Lemma \ref{ColimitSimplicialSet} has the weak homotopy type of the homotopy colimit of the functor $S^f \to \sSet$ given by $I\mapsto NM(I)$. (cf. \cite{Hirschhorn}, Proposition 18.9.4)
\end{Prop}

\begin{pf}
 First, we recall that in order to obtain a model structure on $\Fun(S^f, \sSet)$ where we have a nice description of cofibrations and which satisfies the conditions above, one possibility is to require $S^f$ to be a direct category. This is fulfilled since the assignment $I \mapsto \#I$ gives a linear extension to an ordinal given by, for example, $\# S$. Thus, Theorem 5.1.3 of \cite{Hovey} assures the existence of such a model structure, and it furthermore gives a characterization of cofibrant objects in this model structure. So we only need to check that for each object $I \in S^f$, the induced map $L_I(NM(-)) \to NM(I)$ is a cofibration, where $L_I(F)$ denotes the $I$-latching object of functor $F$. Recall (e.g., from \cite{Hovey}) that $L_I(F)$ is the colimit of the ``restriction'' of $F$ to the category of all non-identity morphisms with target $I$. Note that here, this category is exactly the poset of all proper subsets of $I$. So $L_I(NM(-))$ is by the same argument as in the proof of Lemma \ref{ColimitSimplicialSet} given by
\begin{eqnarray*}
 \bigcup_{J\subsetneq I} NM(J), 
\end{eqnarray*}
and the map is just the inclusion. This yields the claim. 
\end{pf}

Observe that for any $I\subset S$, we have $M(I)\subset M(S)$ and $NM(I)$ is a simplicial subset of $NM(S)$. So we can consider 
\begin{eqnarray*}
 K=\bigcup_{I\in S^f} BM(I) \subset BM(S),
\end{eqnarray*}
realizing the simplicial subset $\widetilde{K}=\bigcup_{I\in S^f} NM(I)$ of $NM$. Our aim is to apply discrete Morse theory to see that the inclusion $K\to BM(S)$ is a homotopy equivalence. The idea is as follows: we will exhibit a proper, acyclic matching on $BM(S)$. This matching will restrict to $K$ and have all essential cells lying in $K$. This will imply that the Morse complex $L$ is the same in both cases. Looking at the situation more closely, we will see that also the Morse equivalence on $K$ is the restriction of the Morse equivalence on $BM(S)$ to $K$. This will imply by two-out-of-three that the inclusion $K\to BM(S)$ is a homotopy equivalence since both $BM(S)\to L$ and its restriction to $K\to L$ are homotopy equivalences. We proceed now by describing the matching on $BM(S)$. 

The Artin monoid $M(S)$ has no non-trivial invertible elements, thus it has Property A by Lemma \ref{MayCriterion}. So we will deal with $Y=r(NM(S))$ instead of the geometric realization, as explained in Remark \ref{RealisierungR}. Recall that the cells of $Y$ are in one-to-one correspondence with the non-degenerate simplices of $NM(S)$, so 
\begin{eqnarray*}
Y^{(*)}=\{[x_n|\ldots | x_1] \, | \,x_i\in M(S), x_i \neq 1 \mbox{ for all } i\}
\end{eqnarray*}
as a set. 

Now we are going to define a matching $\mu_1$ on this set of cells. Recall that we denoted by $\mathcal{D}$ the generating set given by
\begin{eqnarray*}
\mathcal{D}=\{ \Delta_{J}\,|\, \varnothing \neq J\subset S, W(J)\mbox{ is finite} \}. 
\end{eqnarray*}
and for $x\in M$, we denoted $I(x):=\{a\in S|\; \exists \,y\in M\colon \, x=ya\}$. Recall that it is part of Theorem \ref{BSNF} that $\Delta_{I(x)}$ is a right divisor of $x$. 

The rough idea of $\mu_1$ is as follows. Often, to construct a matching on a bar complex, it is helpful to measure to what extent a cell is essential, and this quantity is called ``height'' here. 
The collapsible cells have parts of their boundary as their (redundant) partners, so the latter are determined by an application of a $d_i$ to the former. To obtain an involution, we must be able to recover the original collapsible cell from a redundant one, and we will use the normal form of Theorem \ref{BSNF} to recover the factors from the product. 

We proceed with the precise definition.

\begin{Definition} \label{DefinitionMu1}
We say that an $n$-cell $[x_n|\ldots |x_1]$ of $Y$ is $\mu_1$-\textbf{essential} if for any $1\leq k\leq n$, the product $x_k\ldots x_1$ lies in $\mathcal{D}$. Define $\mu_1([x_n|\ldots |x_1])=[x_n|\ldots |x_1]$ for every essential cell $[x_n|\ldots |x_1]$.

For an arbitrary cell $[x_n|\ldots |x_1]$, we define its $\mu_1$-\textbf{height} by
\begin{eqnarray*}
 \hgta([x_n|\ldots|x_1])=\max\{j|[x_j|\ldots|x_1]\mbox{ is essential}\}
\end{eqnarray*}
If $x_1\notin \mathcal{D}$, set $ \hgta([x_n|\ldots|x_1])=0$. (Set $x_0=1$ for further use.)

For an $n$-cell $[x_n|\ldots |x_1]$ of height $h$ and $1\leq k\leq h$, define $I_k\subset S$ to be the unique subset such that $x_k\ldots x_1=\Delta_{I_k}$. Note that $I_1\subsetneq I_2 \subsetneq \ldots \subsetneq I_h$, and, furthermore, $I_h\subset I(x_{h+1}x_h\ldots x_1)=I(x_{h+1}\Delta_{I_h})$. 

Define an $n$-cell $[x_n|\ldots |x_1]$ of height $h<n$ to be $\mu_1$-\textbf{collapsible} if 
\begin{eqnarray*}
I(x_{h+1}x_h\ldots x_1)=I_h
\end{eqnarray*}
 holds. In this case, we are going to match $[x_n|\ldots|x_1]$ with
 \begin{eqnarray*}
 \mu_1([x_n|\ldots |x_1])=[x_n|\ldots |x_{h+2}|x_{h+1}x_h|x_{h-1}|\ldots |x_1].
\end{eqnarray*}

Define an $n$-cell $[x_n|\ldots |x_1]$ of height $h<n$ to be $\mu_1$-\textbf{redundant} if $I_h\subsetneq J:=I(x_{h+1}x_h\ldots x_1)$. In this case, there exists a unique $y\in M\setminus\{1\}$ such that $x_{h+1}\Delta_{I_h}=y\Delta_J$; furthermore, there is a unique $z\in M\setminus \{1\}$ such that $\Delta_J=z\Delta_{I_h}$. Define 
\begin{eqnarray*}
 \mu_1([x_n|\ldots |x_1])=[x_n|\ldots |x_{h+2}|y|z|x_h|x_{h-1}|\ldots |x_1].
\end{eqnarray*}

In particular, if $[x_n|\ldots |x_1]$ is of height $0$, i.e., if $x_1\notin \mathcal{D}$, we want to define this cell to be redundant (according to our convention, $I(x_0)=\varnothing \subsetneq I(x_1)$). There is then a unique $y\in M\setminus\{1\}$ such that $x_1=y\Delta_{I(x_1)}$; in this case, $z=\Delta_{I(x_1)}$. We define 
\begin{eqnarray*}
 \mu_1([x_n|\ldots |x_1])=[x_n|\ldots |x_{2}|y|\Delta_{I(x_1)}].
\end{eqnarray*}

\end{Definition}

We are now going to show that the assignment $\mu_1$ defines a proper, acyclic matching on $Y$. First, we will show that $\mu_1$ is an involution. This implies that a set of pairs of cells of type 
\begin{eqnarray*}
 (\mu_1-\mbox{ redundant cell}, \mbox{ its collapsible partner})
\end{eqnarray*}
is a candidate for a matching in the sense of Definition \ref{BatziesAcyclicMatching}, and justifies the choice of the names. 

\begin{Lemma}
The map $\mu_1\colon Y^{(*)} \to Y^{(*)}$ from Definition \ref{DefinitionMu1} is an involution.
\end{Lemma}

\begin{pf}
We begin with a collapsible $n$-cell $\underline{x}=[x_n|\ldots |x_1]$ of height $h$. Then 
\begin{eqnarray*}
 \mu_1([x_n|\ldots |x_1])=[x_n|\ldots |x_{h+2}|x_{h+1}x_h|x_{h-1}|\ldots |x_1]
\end{eqnarray*}
 is of height $h-1$, since $x_k\ldots x_1=\Delta_{I_k}$ for $1\leq k\leq h-1$ and $(x_{h+1}x_h)(x_{h-1}\ldots x_1)\notin \mathcal{D}$ by definition. Now since $\underline{x}$ was collapsible, we know that 
 \begin{eqnarray*}
 I((x_{h+1}x_h)(x_{h-1}\ldots x_1))=I_h \supsetneq I_{h-1}.
 \end{eqnarray*}
  Hence, $\mu_1(\underline{x})$ is a redundant cell of height $h-1$. Furthermore, we have $\Delta_{I_h}=x_h\Delta_{I_{h-1}}$ and $x_{h+1}x_h\Delta_{I_{h-1}}=x_{h+1}\Delta_{I_h}$. 

This implies $\mu_1^2(\underline{x})=[x_n|\ldots |x_{h+2}|x_{h+1}|x_h|x_{h-1}|\ldots |x_1]$.

Next, consider a redundant $n$-cell $\underline{x}=[x_n|\ldots |x_1]$ of height $h>0$ with 
\begin{eqnarray*}
 \mu_1([x_n|\ldots |x_1])=[x_n|\ldots |x_{h+2}|y|z|x_h|x_{h-1}|\ldots |x_1].
\end{eqnarray*}
 Then we know that $x_k\ldots x_1=\Delta_{I_k}$ for $1\leq k\leq h$. Furthermore, by definition, $zx_h\ldots x_1=z\Delta_{I_h}=\Delta_{J}$ and $yzx_h\ldots x_1=x_{h+1}x_h\ldots x_1\notin \mathcal{D}$. Thus, the cell $\mu_1(\underline{x})$ has height $h+1$. Moreover, we have $I(yzx_h\ldots x_1)=I(x_{h+1}x_h\ldots x_1)=J$, so that $\mu_1(\underline{x})$ is collapsible. 

Last, we consider the case of a redundant $n$-cell $\underline{x}=[x_n|\ldots |x_1]$ of height $h=0$. We know that $\mu_1([x_n|\ldots |x_1])=[x_n|\ldots |x_{2}|y|\Delta_{I(x_1)}]$. The height of this new cell is at least $1$ since $\Delta_{I(x_1)}\in \mathcal{D}$. Moreover, it is of height exactly $1$ since 
\begin{eqnarray*}
 I(y\Delta_{I(x_1)})=I(x_1)=I(\Delta_{I(x_1)}),
\end{eqnarray*}
and since by assumption $y\Delta_{I(x_1)}=x_1\notin \mathcal{D}$. This implies that the cell $\mu_1([x_n|\ldots |x_1])$ is collapsible of height $1$, and it is mapped by $\mu_1$ to 
\begin{eqnarray*}
 \mu_1^2([x_n|\ldots |x_1])=\mu_1([x_n|\ldots |x_{2}|y|\Delta_{I(x_1)}])=[x_n|\ldots |x_1].
\end{eqnarray*}

We conclude that $\mu_1^2(\underline{x})=[x_n|\ldots |x_{h+2}|x_{h+1}|x_h|x_{h-1}|\ldots |x_1]$. This shows that $\mu_1$ is indeed an involution.
\end{pf}

In order to show that $\mu_1$ defines a matching on $Y^{(*)}$, we still need to show that each $\mu_1$-redundant cell is a regular face of its $\mu_1$-collapsible partner. For this, we are going to exploit the regularity criterion \ref{RegularityCriterion}. According to it, we only need to show that if $\underline{x}=[x_n|\ldots |x_1]$ is a $\mu_1$-collapsible cell of height $h$, then $d_j(\underline{x})\neq d_h(\underline{x})$ for all $0\leq j\neq h\leq n$ (since $d_h(\underline{x})$ is by definition the $\mu_1$-redundant partner of $\underline{x}$). Observe that $h<n$ since $\underline{x}$ is not $\mu_1$-essential, and $h>0$, since all cells of height $0$ are $\mu_1$-redundant.

We have to distinguish several cases. First, assume $1\leq j\neq h \leq n-1$. Without loss of generality, let $j<h$, the other case is treated symmetrically. Then 
\begin{eqnarray*}
d_j([x_n|\ldots |x_1])=[x_n|\ldots |x_{j+2} | x_{j+1}x_j| x_{j-1}| \ldots | x_1].
\end{eqnarray*}
If this term is equal to $d_{h}(\underline{x})$, this implies in particular $x_{j}=x_{j+1}x_j$, since $j<h$. This is a contradiction since $M(S)$ is cancellative and $x_{j+1}\neq 1$ in the non-degenerate simplex $\underline{x}$. So we have to treat the cases $j=0$ and $j=n$. In these cases, $d_j(\underline{x})=d_h(\underline{x})$ would imply $x_{h+1}=x_{h+1}x_h$ or $x_{h}=x_{h+1}x_h$, respectively. This is a contradiction in the same fashion as before. So, by Lemma \ref{RegularityCriterion}, we have indeed a matching on $Y^{(*)}$. We are now going to prove the following. 

\begin{Prop} \label{Mu1ProperAcyclic}
The matching on $Y^{(*)}$ defined by $\mu_1$ is a proper acyclic matching.
\end{Prop}

\begin{pf}
First, we want to show that the matching above is acyclic. Assume we have a cycle 
\begin{eqnarray*}
a_1, a_2, \ldots, a_m=a_1                   
\end{eqnarray*}
in the graph associated to the matching $\mu_1$ on the vertex set $Y^{(*)}$ as in Definition \ref{BatziesAcyclicMatching}. Without loss of generality, we may assume $a_1$ to be a vertex corresponding to a cell of the smallest dimension among $a_1, \ldots, a_m$. Note that each edge in the graph changes the dimension, moreover, the edges decreasing the dimension by $1$ are exactly the ones not in the matching, and the edges increasing the dimension by $1$ are exactly the inverted edges from the matching. So we know that the dimension of $a_2$ has to be $\dim(a_1)+1$, since it is not smaller than $\dim(a_1)$. Thus, $a_1$ and $a_2$ have to be some matched pair, i.e., $\mu_1(a_1)=a_2$. So the cell corresponding to $a_1$ is $\mu_1$-redundant, $a_2$ is $\mu_1$-collapsible and so any edge starting in $a_2$ decreases the dimension. Hence, $\dim(a_3)=\dim(a_1)$ is the smallest dimension in the cycle, so $a_3$ is different from $a_1$ and has to be redundant by the same argument. Therefore, $a_3$ is a redundant boundary of the collapsible partner of the redundant cell $a_1$. Inductively, we obtain a chain 
\begin{eqnarray*}
 a_1\vdash a_3\vdash a_5 \vdash \ldots \vdash a_{2\lfloor \frac{m}{2}\rfloor-1 } \vdash \ldots
\end{eqnarray*}
where $\vdash$ is defined to be the relation for redundant cells $x,z$ of $Y^{(*)}$ with
\begin{eqnarray*}
x \vdash z \Leftrightarrow  z \mbox{ occurs in the boundary of the }\mu_1\mbox{-collapsible partner of } x.
\end{eqnarray*}
So it is enough to show that this relation is noetherian, i.e., has no infinite descending chains. 

 Suppose we have an infinite chain of redundant cells $\underline{x}_1, \underline{x}_2, \ldots$ such that $\underline{x}_{i+1}$ is $d_{k_i}(\mu_1(\underline{x}_i))$ for some $k_i$, and we may assume that $k_i\neq h+1$ so that $\underline{x}_{i+1}\neq \underline{x}_{i}$. Then only finitely many $k_i$ can be $0$ or $n$ since $d_0$ and $d_n$ strictly lower the $S$-length of the product of the entries in the cell label. So we can directly assume there are only $k_i\in \{1,2,\ldots, n-1\}$. We look for possible successors of a redundant $\underline{x}=[x_n|\ldots |x_1]$ of height $h$. For $h+3\leq k\leq n$, the cell $d_k([x_n|\ldots |x_{h+2}|y|z|x_h|x_{h-1}|\ldots |x_1])$ is obviously collapsible, as well as for $1\leq k\leq h$. The cell 
 \begin{eqnarray*}
 d_{h+2}([x_n|\ldots |x_{h+2}|y|z|x_h|x_{h-1}|\ldots |x_1])
 \end{eqnarray*}
  may or may not be redundant, so it is the only possible successor. In any case, note that $\hgta([x_n|\ldots |x_{h+2}y|z|x_h|x_{h-1}|\ldots |x_1])\geq h+1$. Thus, the height in such a chain must strictly increase, so the sequence of redundant cells as above must stabilize after finitely many steps. So we have proven the acyclicity of the matching.
 
 Last, we have to show that the properness of the matching. For this, we want to exploit Lemma \ref{LemmaCompact}. Consider the map 
 \begin{eqnarray*}
 \psi\colon Y^{(*)} & \to & \mathbb{N}\\
{} [x_n|\ldots|x_1] & \mapsto & N_S(x_n\ldots x_1)
\end{eqnarray*}

First, we observe that this is a map of posets: Taking boundaries either leaves the value of $\psi$ constant (if it is $d_i$ for $1\leq i\leq n-1$) or decreases the value (for $i\in \{0,n\}$). Moreover, by definition of $\mu_1$, the value of $\psi$ is the same on the elements matched by $\mu_1$. Last, there are only finitely many elements of $Y^{(*)}$ such that the norm of the product over all entries does not exceed a given value. Thus, by Lemma \ref{LemmaCompact}, the matching is proper. This finishes the proof of the lemma. 
\end{pf}

Next, we observe that all essential cells lie in $K^{(*)}$, a subposet of $Y^{(*)}=BM^{(*)}$. Now if a $\mu_1$-collapsible cell $[x_n|\ldots |x_1]$ lies in some $BM(I)^{(*)}$, so does $d_h([x_n|\ldots |x_1])$, its redundant partner. On the other hand, if $[x_n|\ldots |x_1]$ is $\mu_1$-redundant and lies in $BM(I)^{(*)}$, it is a consequence of the fact that the relations do not change the set of letters of a word that $x_{h+1}=yz$ and $x_{h+1}\in M(I)$ implies $y,z\in M(I)$. So the matching restricts to the subcomplex $K$, and it automatically satisfies the conditions of Definition \ref{BatziesAcyclicMatching} as well as the compactness condition. 

Next, we show that the associated Morse complexes of the matching given by $\mu_1$ and of its restriction to $K$ are the same, and the projections defined in Definition \ref{MorseComplexofMatching} coincide on $K$. Observe that the cells of both Morse complexes $(BM)_{\mu_1}$ and $K_{\mu_1}$ are in one-to-one correspondence with essential cells of either complex, which coincide. Furthermore, it follows inductively from Definition \ref{MorseComplexofMatching} that the projections to the Morse complex coincide on $K$, and this in turn implies that the attaching maps for the Morse complexes coincide. (Here, we also exploit the fact that $K$ is a subcomplex.) Thus we obtain

\begin{Prop} \label{DMTTeilDobrinskaya}
 The inclusion 
\begin{eqnarray*}
 K=\bigcup_{I\in S^f} BM(I) \hookrightarrow BM(S)
\end{eqnarray*}
is a homotopy equivalence. 
\end{Prop}

We are now ready to prove Proposition \ref{MainIngridientDobr}, the missing step in the proof of Dobrinskaya's Theorem \ref{DobrinskayaThm}.
\begin{pf}(of Proposition \ref{MainIngridientDobr})
We put together all the steps done so far. In Proposition \ref{DMTTeilDobrinskaya}, we have shown that $K\simeq BM(S)$. Going through the definition, we observe that $K$ is the geometric realization of the simplicial set $\wt{K}=\bigcup_{I\in S^f} NM(I)$. By Proposition \ref{ModellKatHocolim}, we obtain a weak homotopy equivalence in simplicial sets between $\wt{K}$ and $\hocolim_{S^{f}} NM(-)$. By Theorem \ref{ThomasonTheorem}, this last simplicial set is weakly homotopy equivalent to $N(S^f\int M(-))$, and by Corollary \ref{GrothendieckVergleich}, this simplicial set is in turn weakly homotopy equivalent to $N(S^f\int G(-))$. After geometric realization, we obtain a true homotopy equivalence $K\simeq \abs{N(S^f\int G(-))}$. This completes the proof.
\end{pf}

\begin{remark}
 There are already several applications of discrete Morse theory to hyperplane arrangements in the literature, e.g. in \cite{SalvettiSettepanella}, \cite{MoriSalvetti}, \cite{Delucchi}. Recall that the original formulation of $K(\pi,1)$-conjecture claims that a certain hyperplane arrangement analogue is a $K(G(S),1)$ for an Artin group $G(S)$. 
\end{remark}

\section{Discrete Morse Theory for Chain Complexes}
\label{Discrete Morse Theory for Chain Complexes}

We want to present some further applications of discrete Morse theory to Artin monoids. Unfortunately, the author was so far not able to perform it on the level of topological spaces for the matching of the Section \ref{Squier Complex for Artin Monoids}. 
We introduce an algebraic version here and follow closely the exposition of A.~He\ss{} (\cite{AlexThesis}, see also, e.g., \cite{Kozlov}).

	A \textbf{based chain complex} is a non-negatively graded chain complex $( {C}_*, \partial)$, where each $ {C}_n$ is a free $\IZ$-module, together with a choice of basis $\Omega_n$ for each $ {C}_n$.
	In what follows, $( {C}_*, \Omega_*, \partial)$ will always be a based chain complex.

	We equip each $ {C}_n$ with the inner product $\langle \variable , \variable \rangle\colon  {C}_n \times  {C}_n \to \IZ$ obtained by regarding $\Omega_n$ as an orthonormal basis for ${C}_n$.
	If $x$, $y$ have the ``wrong'' dimensions, i.e., if $x\in C_n$, but $y\notin C_{n}$, then we set their product $\langle  x, y \rangle$ to be zero.

	\begin{Definition}
		A $\mathbb{Z}$-compatible \textbf{matching} on a based chain complex $( {C}_*, \Omega_*, \partial)$ is an involution $\mu\colon \Omega_* \to \Omega_*$ satisfying the following property: For every $x \in \Omega_*$ which is not a fixed point of $\mu$, we have $\langle \partial x, \mu(x) \rangle = \pm 1$ or $\langle \partial \mu(x), x \rangle = \pm 1$. (This last condition is called $\mathbb{Z}$-compatibility.)
		
		The fixed points of a matching $\mu\colon \Omega_* \to \Omega_*$ are called \textbf{essential}. If $x \in \Omega_n$ is not a fixed point, then $\mu(x) \in \Omega_{n-1} \cup \Omega_{n+1}$. We say that $x$ is \textbf{collapsible} if $\mu(x) \in \Omega_{n-1}$, and it is called \textbf{redundant} if $\mu(x) \in \Omega_{n+1}$.
	\end{Definition}

	\begin{remark}\label{Zcompatibility}
Let $\mu\colon \Omega_* \to \Omega_*$ be an involution. Assume we know that all non-fixed points of $\mu$ are either collapsible or redundant. Then it is enough to check $\Mf{\partial\mu(x), x}=\pm 1$ for redundant cells in order to check that  $\mu$ is $\mathbb{Z}$-compatible. Indeed, let $x \in \Omega_n$ be a non-fixed point of an involution $\mu$ as above. We have to show that $\Mf{\mu(x), \partial x}=\pm 1$ for the case that $x$ is collapsible. In this case, the image $\mu(x)\in \Omega_{n-1}$ is redundant since $\mu(\mu(x))=x$ is in $\Omega_n$. So we know that for $y=\mu(x)$, we have $\Mf{\partial\mu(y), y}=\pm 1$. Inserting $y=\mu(x)$, we obtain $\Mf{\mu(x), \partial x}=\pm 1$. 
	\end{remark}

	Let $\mu$ be a matching on $( {C}_*, \Omega_*, \partial)$. For two redundant basis elements $x$, $z \in \Omega_*$ set $x \vdash z$ to be the relation ``$z$ occurs in the boundary of the collapsible partner of $x$'', i.e. $\Mf{\partial \mu(x), z} \neq 0$.

	\begin{Definition}\label{AlgebraicMatchingNoetherian}
		A matching on a based chain complex is called \textbf{noetherian} if every infinite chain $x_1 \vdash x_2 \vdash x_3 \vdash \ldots$ eventually stabilizes.
	\end{Definition}
	
	\begin{remark}
	This definition of noetherianity for a matching on based chain complexes is not the standard one. Yet, as observed in \cite{AlexThesis}, §1.1, this definition is equivalent to the usual one and is often easier to check. 
	\end{remark}
	
	Given a noetherian matching $\mu$ on $( {C}_*, \Omega_*, \partial)$, we define a linear map $\theta^\infty\colon  {C}_* \to  {C}_*$ as follows. Let $x \in \Omega_*$. If $x$ is essential, we set $\theta(x) = x$. If $x$ is collapsible, we set $\theta(x) = 0$, and if $x$ is redundant we set $\theta(x) = x - \varepsilon \cdot \partial\mu(x)$, where $\varepsilon = \Mf{\partial \mu(x),x}$.
	
	Note that, if $x$ is redundant, then $\langle x, \theta(x) \rangle = 0$. It is now not hard to check that for every $x \in \Omega_*$ the sequence $\theta(x), \theta^2(x), \theta^3(x), \ldots$ stabilizes (cf. also \cite{AlexThesis}, Section 1.1), and we define $\theta^\infty(x) := \theta^N(x)$ for $N$ large enough. We linearly extend this map to obtain $\theta^\infty\colon  {C}_* \to  {C}_*$.
	
	We can now state the main theorem of discrete Morse theory for chain complexes.
	
	\begin{Theorem}[Brown, Cohen, Forman]\label{BCF}
		Let $( {C}_*, \Omega_*, \partial)$ be a based chain complex and let $\mu$ be a noetherian matching on it. Denote by $ {C}_*^\theta = \im(\theta^\infty\colon  {C}_* \to  {C}_*)$ the $\theta$-invariant chains. Then $( {C}_*^\theta, \theta^\infty \circ \partial|_{\im(\theta^{\infty})})$ is a chain complex, and the map
		\begin{align*}
			\theta^\infty\colon ( {C}_*, \partial) \longrightarrow \left( {C}_*^\theta, \theta^\infty \circ \partial|_{\im(\theta^{\infty})}\right)
		\end{align*}
		is a chain homotopy equivalence. A basis of $ {C}_*^{\theta}$ is given by the essential cells. 
	\end{Theorem}
	
	For a proof see e.g. \cite{Forman}.

We want to pin down the connection between both flavors of discrete Morse theory introduced so far. 

\begin{Lemma}\label{DMTtoptoalg}
Let $X$ be a CW complex and let $A$ be a proper, acyclic matching on $X^{(*)}$. Then it induces a noetherian, $\mathbb{Z}$-compatible matching on the (based) cellular chain complex $C_*(X)$, and the essential cells of both matchings coincide.
\end{Lemma}

\begin{pf}
Recall that the cellular chain complex $C_*(X)$ has exactly $X^{(*)}$ as a basis. We define the map $\mu\colon X^{(*)}\to X^{(*)}$ by
\begin{eqnarray*}
\mu(x)=\begin{cases}
x, \mbox{ if } x \mbox{ is } A-\mbox{essential},\\
y, \mbox{ if } (x\to y)\in A \mbox{ or } (y\to x) \in A.
\end{cases}
\end{eqnarray*}
This is obviously a well-defined involution, and the notions $A$-essential, $A$-collapsible, $A$-redundant and $\mu$-essential, $\mu$-collapsible, $\mu$-redundant coincide. 

Next, we have to check that $\mu$ is $\mathbb{Z}$-compatible. As explained in Remark \ref{Zcompatibility}, we only have to show $\Mf{\partial \mu(x),x}=\pm 1$ for redundant $n$-cells $x$. By definition, the cell $\mu(x)$ is $n+1$-dimensional, and $x$ is a regular face of $\mu(x)$. Now $\Mf{\partial \mu(x),x}$ is given in the cellular chain complex by the degree of the map 
\begin{eqnarray*}
S^n \xrightarrow{f_{\mu(x)}|_{S^n}} X^{(n)} \to S^n,
\end{eqnarray*}
where $f_{\mu(x)}\colon D^{n+1} \to X$ is the characteristic map of the cell $\mu(x)$, so that $f_{\mu(x)}|_{S^n}$ is the attaching map of $\mu(x)$, and the second map is given by collapsing everything in the $n$-skeleton $X^{(n)}$ outside the open cell $\mathring{x}$ to a point. Now since $x$ is a regular face of $\mu(x)$, we know that the pre-image $f_{\mu(x)}^{-1}(\mathring{x})$ is mapped homeomorphically to $\mathring{x}$. Using local degree calculation, we may conclude that $\Mf{\partial \mu(x),x}=\pm 1$, as desired. 

Last, we want to show that the matching $\mu$ is noetherian. So assume we had an infinite chain  $x_1 \vdash x_2 \vdash x_3 \vdash \ldots$. Observe that whenever $x \vdash y$, then there is a path from $x$ to $y$. So, every $x_i$ lies in $X^{(*)}_{\preceq x_1}$, where by $x_1$, we mean here the element of $A^{(*)}$ containing $x_1$. So by the properness of the matching, there can be only finitely many different $x_i$. Since the chain is infinite, it has to contain a cycle, which would exactly correspond to a cycle in the associated graph $G_{X}^{A}$. This is a contradiction, since the matching $A$ was assumed to be acyclic. This completes the proof of the noetherianity and thus of the lemma.
\end{pf}

\section{Squier Complex for Artin Monoids}
\label{Squier Complex for Artin Monoids}

In this section, we consider first the noetherian matching $\mu_1$ on the bar complex of an Artin monoid, induced by the matching $\mu_1$ on $BM(S)$. We construct a further noetherian matching $\mu_2$ on the obtained chain complex, so that the resulting chain complex is related to the one defined by Squier (\cite{Squier}). 

Let $M:=M(S)$ be again an Artin monoid with Artin-Coxeter generating set $S$. Let $W(S)$ be the corresponding Coxeter group. Set again
\begin{eqnarray*}
 S^{f}=\{I\subset S| W(I)\mbox{ is finite} \}.
\end{eqnarray*}
Furthermore, set again $\mathcal{D}=\{\Delta_{I}:=\llcm(I)| \varnothing \neq I\in S^{f}\}$.  

\begin{Lemma}
There is a noetherian matching $\mu_1$ on the bar complex of the Artin monoid $M(S)$ with the same essential, collapsible and redundant cells as the ones described in Definition \ref{DefinitionMu1}.
\end{Lemma}

\begin{pf}
This follows directly from Proposition \ref{Mu1ProperAcyclic} and Lemma \ref{DMTtoptoalg}. 
\end{pf}

\begin{remark}\label{DifferentialSquier1}
By Theorem \ref{BCF}, we know that the complex $(C_*^{\theta_1}, d_*^{\theta_1})$ computes the homology of an Artin monoid $M$, where $C_n^{\theta_1}$ has as a $\mathbb{Z}$-basis the $\mu_1$-essential $n$-cells, and $d_*^{\theta_1}=\theta_1^{\infty}\circ d$. Now if $\underline{x}=[x_n|\ldots |x_1]$ is an essential cell, it is clear that $d_i(\underline{x})$ is essential for $1\leq i\leq n$, while $d_0(\underline{x})$ may or may not be essential. Thus, we have 
\begin{eqnarray*}
d_*^{\theta_1}([x_n|\ldots |x_1])=d(\underline{x})-[x_n|\ldots |x_2]+\theta_1^{\infty}([x_n|\ldots |x_2])
\end{eqnarray*}
Note that the summands of $d_*^{\theta_1}(\underline{x})$ are either $\pm d_i(\underline{x})$ for $1\leq i\leq n-1$ or the product of their entries have smaller $S$-length than such of $\underline{x}$. We will need this description later.

Furthermore, note that any $\mu_1$-essential cell $[x_n|\ldots |x_1]$ is uniquely characterized by the sequence $I_n\supsetneq I_{n-1} \supsetneq \ldots \supsetneq I_1$. We denote the set of such cells in dimension $n$ by $\Omega_n^{\theta_1}$. (This is a basis for $C_n^{\theta_1}$.)
\end{remark}

Now we are going to define a noetherian matching $\mu_2$ on the obtained chain complex $(C_*^{\theta_1}, d_*^{\theta_1})$ making it smaller again. For this, choose any linear order $<$ on the set $S$. 

The rough idea of $\mu_2$ is as follows. We want the essential cells of $\mu_2$ to be those where the sets in the characterizing sequence grow one element at a time, and in addition, the new element is assumed to be larger then the old ones. Again, we define a height function measuring up to which point the beginning of the characterizing sequence is already essential. Then, we call all the cells redundant where we can enlarge this starting sequence by borrowing the maximal element from the next set in the characterizing sequence. If this is not possible, the cell is collapsible and we will forget about the last set in the starting essential sequence.

We now describe the essential, collapsible and redundant cells of the matching $\mu_2$.

\begin{Definition}
Let $[x_n|\ldots |x_1]$ be an $n$-cell in $\Omega_n^{\theta_1}$. We say $[x_n|\ldots |x_1]$ to be $\mu_2$-\textbf{essential} if for any $1\leq k\leq n$, $I_{k}\setminus I_{k-1}=\{a_k\}$ and $a_k=\max I_k$. (Here, we set $I_0=\varnothing$.) Define $\mu_2([x_n|\ldots |x_1])=[x_n|\ldots |x_1]$ for an essential cell $[x_n|\ldots |x_1]$.

For an arbitrary cell $[x_n|\ldots |x_1]$, we define its $\mu_2$-\textbf{height} by
\begin{eqnarray*}
 \hgtb([x_n|\ldots|x_1])=\max\{j|[x_j|\ldots|x_1]\mbox{ is essential}\}
\end{eqnarray*}
If $\# I_1 > 1$, set $ \hgtb([x_n|\ldots|x_1])=0$.

Define an $n$-cell $[x_n|\ldots |x_1]$ of height $h<n$ to be $\mu_2$-\textbf{collapsible} if $\max I_{h+1} =\max I_h$ holds. In this case, set \begin{eqnarray*}
 \mu_2([x_n|\ldots |x_1])=[x_n|\ldots |x_{h+2}|x_{h+1}x_h|x_{h-1}|\ldots |x_1].
\end{eqnarray*}
The characterizing sequence of the new element is 
\begin{eqnarray*}
I_n\supsetneq \ldots \supsetneq I_{h+2}\supsetneq I_{h+1}\supsetneq I_{h-1}\supsetneq \ldots \supsetneq I_{1}.
\end{eqnarray*}

Define an $n$-cell $[x_n|\ldots |x_1]$ of height $h<n$ to be $\mu_2$-\textbf{redundant} if 
\begin{eqnarray*}
b:=\max I_{h+1} > a_h=\max I_h.
\end{eqnarray*}
 Observe that in this case $\# I_{h+1}\geq 2+\# I_h$ since otherwise the cell would have at least height $h+1$. Thus there exist $u,v\in M\setminus\{1\}$ such that $\Delta_{I_{h+1}}=u\Delta_{I_{h}\cup \{b\}}$ and $\Delta_{I_{h}\cup \{b\}}=v\Delta_{I_h}$. Define 
\begin{eqnarray*}
 \mu_2([x_n|\ldots |x_1])=[x_n|\ldots |x_{h+2}|u|v|x_h|x_{h-1}|\ldots |x_1].
\end{eqnarray*}
Note that the characterizing sequence of the new element is 
\begin{eqnarray*}
I_n\supsetneq \ldots \supsetneq I_{h+2}\supsetneq I_{h+1}\supsetneq I_{h}\cup \{b\} \supsetneq I_h\supsetneq I_{h-1}\supsetneq \ldots \supsetneq I_{1}.
\end{eqnarray*}
Observe furthermore that $\Delta_{I_{h+1}}=u\Delta_{I_{h}\cup \{b\}}=uv\Delta_{I_h}$ implies $x_{h+1}=uv$. 
\end{Definition}

We are going to prove that $\mu_2$ is a noetherian matching on $(C_*^{\theta_1}, d_*^{\theta_1})$.  

\begin{Prop} \label{MatchingMu2}
 For any Artin monoid $M$, the map $\mu_2\colon \Omega_*^{\theta_1}\to \Omega_*^{\theta_1}$ defined as above gives a noetherian, $\mathbb{Z}$-compatible matching on $(C_*^{\theta_1}, d_*^{\theta_1})$. 
\end{Prop}

\begin{pf}
First, we are going to show that $\mu_2$ is an involution. We begin with a collapsible $n$-cell $\underline{x}=[x_n|\ldots |x_1]$ of height $h$. Then 
\begin{eqnarray*}
 \mu_2([x_n|\ldots |x_1])=[x_n|\ldots |x_{h+2}|x_{h+1}x_h|x_{h-1}|\ldots |x_1] 
\end{eqnarray*}
is of height $h-1$, since $I_{k}\setminus I_{k-1}=\{a_k\}$ and $a_k=\max I_k$ for  $1\leq k\leq h-1$ and $\# (I_{h+1}\setminus I_{h-1})\geq 2$. Since $\underline{x}$ was collapsible of height $h$, we know that $\max I_{h+1}=a_h>\max I_{h-1}=a_{h-1}$, so that $\mu_2(\underline{x})$ is a redundant cell of height $h-1$. Furthermore, we have $I_h=I_{h-1}\cup\{a_h\}$ and $\Delta_{I_h}=x_h\Delta_{I_{h-1}}$ as well as $\Delta_{I_{h+1}}=x_{h+1}\Delta_{I_h}$. This implies by definition $\mu_2([x_n|\ldots |x_1])=[x_n|\ldots |x_{h+2}|x_{h+1}|x_h|x_{h-1}|\ldots |x_1]$.

Next, consider a redundant $n$-cell $\underline{x}=[x_n|\ldots |x_1]$ of height $h$ with 
\begin{eqnarray*}
\mu_2([x_n|\ldots |x_1])=[x_n|\ldots |x_{h+2}|u|v|x_h|x_{h-1}|\ldots |x_1],
\end{eqnarray*}
 as defined above. Again, we have $I_{k}\setminus I_{k-1}=\{a_k\}$ and $a_k=\max I_k$ for  $1\leq k\leq h$. Furthermore, we know that $vx_h\ldots x_1=\Delta_{I_{h}\cup \{b\}}$ and $b=\max(I_{h}\cup \{b\}) >\max I_h$ by definition. In addition, we have that $b=\max(I_{h}\cup \{b\})=\max I_{h+1}$, so that $\mu_2(\underline{x})$ is a collapsible cell of height $h+1$. 

We conclude that $\mu_2^2(\underline{x})=[x_n|\ldots |x_{h+2}|x_{h+1}|x_h|x_{h-1}|\ldots |x_1]$. This shows that $\mu_2$ is indeed an involution. 

Now we observe that $\mu_2$ is $\mathbb{Z}$-compatible. Consider a redundant cell $\underline{x}:=[x_n|\ldots |x_1]$. Clearly, $d_{h+1}([x_n|\ldots |x_{h+2}|u|v|x_h|x_{h-1}|\ldots |x_1])=\underline{x}$ if $\hgtb(\underline{x})=h$ and it is easy to see that none of the other $d_i$-summands for $1\leq i \leq n$ produces $\underline{x}$. Since the other summands of $d^{\theta_1}_*$ have smaller $S$-norm, they cannot coincide with $\underline{x}$. This shows that $\mu_2$ is a $\mathbb{Z}$-compatible matching.

Finally, we are going to show that the matching $\mu_2$ is noetherian. Suppose we have an infinite sequence of redundant $n$-cells $\underline{x}_1, \underline{x}_2, \ldots$, such that $\underline{x}_{i+1}$ is a summand of $d^{\theta_1}_*(\mu_2(\underline{x}_i))$. We may assume that $\underline{x}_{i+1}\neq \underline{x}_{i}$. Moreover, we may assume that the $S$-length of the product of all entries is constant, since it is non-increasing and finite. Thus, we may assume that $\underline{x}_{i+1}=d_{k_i}(\mu_2(\underline{x}_i))$ with $k_i\in \{1,2,\ldots, n\}$.

Define a (lexicographic) order on characterizing sequences as follows: Let 
\begin{eqnarray*}
 (A_m\supsetneq A_{m-1} \supsetneq \ldots \supsetneq A_1)<(B_l\supsetneq B_{l-1} \supsetneq \ldots \supsetneq B_1)
\end{eqnarray*}
hold if there exist $s<\max\{m,l\}$ such that $\max A_i=\max B_i$ for all $1\leq i\leq s-1$ and $\max A_s <\max B_s$, or if $\max A_i=\max B_i$ for all $1\leq i\leq \max\{m,l\}$ and there is a $t<\max\{m,l\}$ such that $\# A_i=\# B_i$ for $1\leq i\leq t-1$ and $\# A_t> \# B_t$. (For $i>m$, set $A_i=A_m$, similar for $B_i$. The opposite directions of the inequality signs may be a bit confusing at the beginning, but this is the definition we will need.)

 We consider the characterizing sequences for the successors of a redundant cell $\underline{x}=[x_n|\ldots |x_1]$ of height $h$. For $h+3\leq k\leq n$, the cell 
 \begin{eqnarray*}
 d_k([x_n|\ldots |x_{h+2}|u|v|x_h|x_{h-1}|\ldots |x_1])
 \end{eqnarray*}
  is obviously collapsible and thus not a successor of $\underline{x}$. For $1\leq k\leq h$, the characterizing sequence of $d_k([x_n|\ldots |x_{h+2}|u|v|x_h|x_{h-1}|\ldots |x_1])$ is larger than the one of $\underline{x}$ since both coincide for $1\leq i\leq k-1$ and for $i=k$ we have $\max I_{k+1}=a_{k+1}>\max I_k=a_k$. (This holds also for $k=h$ since also $b=\max(I_{h}\cup \{b\}) >\max I_h$.)
 The cell $d_{h+2}([x_n|\ldots |x_{h+2}|u|v|x_h|x_{h-1}|\ldots |x_1])=[x_n|\ldots |x_{h+2}u|v|x_h|x_{h-1}|\ldots |x_1]$ has the characterizing sequence 
 \begin{eqnarray*}
I_n\supsetneq \ldots \supsetneq I_{h+2}\supsetneq I_{h}\cup \{b\} \supsetneq I_h\supsetneq I_{h-1}\supsetneq \ldots \supsetneq I_{1}
\end{eqnarray*}
This sequence has in each place the same maximum as the original one and 
\begin{eqnarray*}
\#(I_{h}\cup \{b\})<\#I_{h+1}
\end{eqnarray*} 
holds, so that this sequence is again larger than the one of $\underline{x}$. 

 Thus, the characterizing sequence in such a chain must strictly increase, so the sequence of redundant cells as above must stabilize after finitely many steps.
\end{pf}

Note that the essential cells in dimension $k$ of the new complex are in one-to-one correspondence with the $k$-element subsets of $S$ which lie in $S^f$. This reproves the following proposition found by Squier (\cite{Squier}):

\begin{Prop}(\cite{Squier}, Theorem 7.5)
 Let $M(S)$ be an Artin monoid. Then there is a complex which computes the homology of $M(S)$, consisting of free modules $C_k$ with a basis given by the $k$-element subsets of $S$ which lie in $S^f$. 
\end{Prop}

\begin{remark}
\begin{enumerate}
 \item The complex given by Squier in \cite{Squier} has the advantage that the differentials can be described explicitly. Assume $S$ carries a linear order $<$. The value of the differential on the generator $[I]=[a_1<a_2<\ldots <a_k]$ is given by
\begin{eqnarray*}
 \partial([I])=\sum_{i=1}^{k} (-1)^{i-1}\cdot \left(\sum_{uv=\Delta_{I}\Delta_{I\setminus\{a_i\}}^{-1}} (-1)^{N_S(u)}\right) \left[I\setminus\{a_i\}\right].
\end{eqnarray*}
We will call this complex the \emph{Squier complex}. 
\item There is at least one further complex in the literature computing the homology of Artin monoids which has the same number of generators as the Squier complex: It comes from a space homotopy equivalent to $BM$ and is often called the \emph{Salvetti complex}. It can be found e.g. in \cite{ConciniSalvetti}, \cite{Salvetti}, cf. also \cite{CD}. There seems to be no account in the literature comparing the Salvetti complex and the Squier complex  (cf. Introduction of \cite{EllisSkoeldberg}). It is also unclear to the author how the differentials in the complex obtained after applying the matching $\mu_2$ of this section compare to those of the Salvetti complex and to those of the Squier complex. 
\end{enumerate}
\end{remark}

\section{Generalized Charney-Meier-Whittlesey Complex}
\label{Matching on the bar complex}

The results of this section were found partially joint with A. He\ss{}. The aim of this section is to generalize a theorem by R.~Charney, J.~Meier and K.~Whittlesey (\cite{CMW}), which gives a small classifying space model for Artin groups of finite type, to all Artin monoids. Indeed, in \cite{CMW}, the more general setting of Garside groups is considered, and the generalization goes through the same way for left cancellative, locally left Gaussian monoids which are atomic. For the definitions, we refer to \cite{DehornoyLafont}. The details of the general statement can be found in the author's thesis (\cite{MyThesis}). Here, we give a proof for Artin monoids using discrete Morse theory. A proof of the original theorem of \cite{CMW}  via discrete Morse theory can be found in \cite{AlexThesis}.

We will again make use of the divisibility properties in an Artin monoid, and also of several properties as written down or inspired by \cite{DehornoyLafont}. We begin with a definition:

\begin{Definition} \label{DefinitionLeftComplement}
 Let $M$ be a left-cancellative monoid with no non-trivial invertible elements. For any two $x,y\in M$ whose left least common multiple exists, we denote by $x/y$ the unique element such that $\llcm(x, y)=(x/y)\cdot y$, called the \textbf{left complement} of $y$ in $x$. Here, $\llcm$ denotes the left least common multiple of two elements.
\end{Definition}

We recall furthermore the following property of Artin monoids:

\begin{Lemma}[\cite{BrieskornSaito}, Prop. 4.1] \label{ArtinGauss}
Any two elements of an Artin monoid admit either a left least common multiple or do not admit left common multiples at all. 
\end{Lemma}

The proof of the following statement follows closely the proof of Lemma 1.7 of \cite{DehornoyLafont}.

\begin{Lemma}
 Let $M=M(S)$ be an Artin monoid. Let $\mathcal{E}$ be a generating set closed under left complement and left least common multiple, i.e., whenever $p,q\in\mathcal{E}$ and their left least common multiple $y$ exists, then $y, (p/q), (q/p) \in \mathcal{E}\cup\{1\}$. Let $a\neq 1$ be an element of $M$, and let $b$ be in $\mathcal{E}\cup\{1\}$. Then there exists a (unique) greatest right divisor $d$ of $a$ for which $db \in \mathcal{E}$. 
\end{Lemma}

\begin{pf}
 First, we observe that the set 
 \begin{eqnarray*}
 A=\{z\in M\,|\, \mbox{There exists } u\in M \mbox{ with } a=zu\mbox{ and } ub\in \mathcal{E}\cup\{1\}\}
 \end{eqnarray*}
  is non-empty since $b=1\cdot b\in \mathcal{E}$. Consider an element $c\in A$ of minimal $S$-length and write $a=cd$. Note that by definition $db\in \mathcal{E}\cup\{1\}$. We are going to show that $d$ has the desired property. Let $a=uv$ for $u,v\in M$ such that $vb\in \mathcal{E}\cup\{1\}$. We have to show that $v$ is a right divisor of $d$. 

Observe that $vb$ is the left least common multiple of $b$ and $vb$. Thus, $v$ is a left complement of those (cf. Definition \ref{DefinitionLeftComplement}) and has to lie in $\mathcal{E}\cup\{1\}$. Similarly, we observe that $d \in \mathcal{E}\cup\{1\}$. Since $a$ is a common left multiple of $d$ and $v$, they must have a left least common multiple $sd=tv$ by Lemma \ref{ArtinGauss} and $sd, s,t\in \mathcal{E}\cup\{1\}$ by assumption. Moreover, there exists $x\in M$ such that $a=x(sd)$ by the definition of a least common multiple. Furthermore, it is easy to see that $sdb=tvb$ is the left least common multiple of $db$ and $vb$, thus lies again in $\mathcal{E} \cup \{1\}$. So $x\in A$ and, since $a=xsd=cd$, we have $xs=c$ and $x$ is a left divisor of $c$. By the minimality of $c$, we have $s=1$, so $d=tv$; thus any right-divisor $v$ of $a$ so that $vb\in \mathcal{E}\cup\{1\}$ is a right-divisor of $d$. This yields the claim.
\end{pf}

\begin{Notation}
 Let $M$ be an Artin monoid, $\mathcal{E}$ a generating set closed under left complement and left least common multiple. Let $a\neq 1$ be an element of $M$, and let $b$ be in $\mathcal{E}\cup \{1\}$. We write $\gamma(a,b)$ for the greatest right divisor $d$ of $a$ for which $db \in \mathcal{E}$. Furthermore, we write $\psi(a,b)$ for the unique element with $a=\psi(a,b)\gamma(a,b)$.
\end{Notation}

	We are now going to construct a third proper, acyclic matching on the classifying space $BM(S)$ of an Artin monoid $M(S)$. 
   
	Let $\mathcal{E}$ be a generating set for $M=M(S)$, closed under left complement and left least common multiple. Define a subset $\mathcal{E}_n$ of cells of $BM^{(*)}$ by 
\begin{eqnarray*}
 \mathcal{E}_n=\{[x_n|\ldots|x_1]\in BM^{(*)} \left| \mbox{ For all } 1\leq k \leq n, x_k\ldots x_1 \in \mathcal{E}\right\}.
\end{eqnarray*}

	\begin{Prop}\label{matching}
		There exists a proper, acyclic matching $\mu$ on $BM^{(*)}$ with the property that an $n$-cell $x \in BM^{(*)}$ is a $\mu$-essential cell if and only if $x \in \mathcal{E}_n$.
	\end{Prop}

	\begin{pf}
		First, we define the height of a cell $[x_n | \ldots |x_1] \in BM^{(*)}$ to be the maximal integer $h \geq 0$ subject to $[x_h | \ldots | x_1] \in \mathcal{E}_h$. If $h = n$, then $\mu$ is defined to fix this element. Otherwise, $h+1\leq n$ and by definition $x_{h+1} x_h \ldots x_1 \notin \mathcal{E}$. For convenience, set $x_0=1$. We now distinguish two cases.
		\begin{enumerate}
			\item	
				If $\gamma(x_{h+1}, x_{h}\ldots x_1x_0) = 1$, then we call the cell $[x_n|\ldots |x_1]$ collapsible and set
				\begin{align*}
					\mu([x_n | \ldots | x_1]) = [x_n|\ldots|x_{h+2}|x_{h+1}x_{h}|x_{h-1}|\ldots|x_1].
				\end{align*}
			Observe that the new cell has height $h-1$.
			\item	
				If $\gamma(x_{h+1}, x_{h}\ldots x_1x_0) \neq 1$, then we call $[x_n|\ldots |x_1]$ redundant and set 
				\begin{eqnarray*}
					&&a = \psi(x_{h+1}, x_{h}\ldots x_1),\\
					&&d = \gamma(x_{h+1}, x_{h}\ldots x_1) \mbox{ and }\\
					&&\mu([x_n | \ldots | x_1]) = [x_n|\ldots|x_{h+2}|a|d|x_h|\ldots|x_1].
				\end{eqnarray*}
				Note that $a \neq 1$, because by definition $dx_h \ldots x_1 \in \mathcal{E}$ but 
\begin{eqnarray*}
 adx_h \ldots x_1=x_{h+1}x_h\ldots x_1 \notin \mathcal{E}.
\end{eqnarray*}
Furthermore, by this argumentation we see that the new cell has height $h+1$. In particular, if we started with a cell of height $0$, we will get into this case since $\gamma(x_1, x_0)=\gamma(x_1,1)$ is exactly the greatest divisor of $x_1$ lying in $\mathcal{E}$. The element $\gamma(x_1, 1)$ is non-trivial since $x_1\neq 1$. 
		\end{enumerate}
We are now going to show step by step that $\mu$ defines a proper, acyclic matching on $BM^{(*)}$.

Our first goal is to show that $\mu$ is an involution. This will imply that the non-fixed points of $\mu$ give a collection of cell pairs, in each of which the redundant cell is a face of its collapsible partner, and by definition, each cell of $BM^{(*)}$ appears at most once in this collection. 
		
		Let $\underline{x} = [x_n|\ldots|x_1]$ be redundant of height $h$. We will first show that 
		\begin{eqnarray*}
		\mu(\underline{x}) = [x_n|\ldots|x_{h+2}|a|d|x_h|\ldots|x_1]
		\end{eqnarray*}
		 is collapsible of height $h+1$.

		Set $c=\gamma(a, dx_h\ldots x_1)$. Then $a=yc$ for $y=\psi(a, dx_h\ldots x_1)$, so $x_{h+1}=ad=y(cd)$ and $(cd)x_h\ldots x_1 \in \mathcal{E}$. By definition of $d$, we have $cd=d$ and $c=1$.
		Thus, $\mu(\underline{x})$ is collapsible of height $h+1$. Hence,
		\begin{align*}
			\mu^2(\underline{x}) = \mu([x_n|\ldots|x_{h+2}|a|d|x_{h}|\ldots |x_1]) = [x_n|\ldots|x_{h+2}|ad|x_{h}|\ldots |x_1] = \underline{x}.
		\end{align*}

		Now let $\underline{x} = [x_n|\ldots |x_1]$ be collapsible of height $h$. We will first show that 
		\begin{eqnarray*}
		\mu(\underline{x}) = [x_n|\ldots|x_{h+2}|x_{h+1}x_{h}|x_{h-1}|\ldots|x_1]
		\end{eqnarray*}
		 is redundant of height $h-1$. 
		
		We have to compute $u=\gamma(x_{h+1}x_{h}, x_{h-1}\ldots x_1)$. Observe that $x_{h}$ is a right divisor of $x_{h+1}x_h$ and $x_h(x_{h-1}\ldots x_1)\in \mathcal{E}$ by assumption on $\underline{x}$. So by the definition of $u$, we have $u=sx_{h}$ for some $s\in M$, and $x_{h+1}x_h=ru=r(sx_{h})$ for $r =\psi(x_{h+1}x_h, x_{h-1}\ldots x_1)$. Thus, $x_{h+1}=rs$ and $u(x_{h-1}\ldots x_1)=s(x_hx_{h-1}\ldots x_1)\in \mathcal{E}$. By definition, $s$ is a right divisor of $\gamma(x_{h+1}, x_h\ldots x_1)$ which is $1$ since $\underline{x}$ was collapsible. This implies $s=1$ and $\gamma(x_{h+1}x_{h}, x_{h-1}\ldots x_1)=x_h$ since there are no non-trivial invertible elements in $M$.

		Since $x_h \neq 1$, this proves that $\mu(\underline{x})$ is redundant of height $h-1$. Hence,
		\begin{align*}
			\mu^2(\underline{x}) = \mu([x_n|\ldots|x_{h+2}|x_{h+1}x_{h}|x_{h-1}|\ldots|x_1]) = 
			\underline{x}.
		\end{align*}
		This shows that $\mu$ is an involution.
		
		Now we have to check that any redundant cell is a regular face of its collapsible partner. For doing so, we are going to exploit again the regularity criterion \ref{RegularityCriterion}. According to it, we only need to show that if $\underline{x}=[x_n|\ldots |x_1]$ is a collapsible cell of height $h$, then $d_j(\underline{x})\neq d_h(\underline{x})$ for all $0\leq j\neq h\leq n$ (since $d_h(\underline{x})$ is by definition the redundant partner of $\underline{x}$). Observe that $h<n$ since $\underline{x}$ is not essential, and $h>0$, since all cells of height $0$ are redundant.

We have to distinguish several cases, similarly as in Section \ref{A reformulation of the Kpi1-conjecture}. First, assume $1\leq j\neq h \leq n-1$. Without loss of generality, let $j<h$, the other case is treated symmetrically. Then 
\begin{eqnarray*}
d_j([x_n|\ldots |x_1])=[x_n|\ldots |x_{j+2} | x_{j+1}x_j| x_{j-1}| \ldots | x_1].
\end{eqnarray*}
If this term is equal to $d_{h}(\underline{x})$, this implies in particular $x_{j}=x_{j+1}x_j$, since $j<h$. This is a contradiction since $M(S)$ is cancellative and $x_{j+1}\neq 1$ in the non-degenerate simplex $\underline{x}$. So we have to treat the cases $j=0$ and $j=n$. In these cases, $d_j(\underline{x})=d_h(\underline{x})$ would imply $x_{h+1}=x_{h+1}x_h$ or $x_{h}=x_{h+1}x_h$, respectively. This is a contradiction in the same fashion as before. So, by Lemma \ref{RegularityCriterion}, we have indeed a matching on $BM^{(*)}$.  		
		
		Next, we are going to check that the matching $\mu$ is acyclic. As before, it is enough to check that there is no infinite sequence of redundant cells $\underline{x}_1 \vdash \underline{x}_2 \vdash \ldots$. 
		
 Since $d_0$ and $d_n$ strictly decrease the $S$-length of the product of the entries of a cell label, we may assume that each $\underline{x}_i$ in such a sequence is given by $d_{k_i}(\mu(\underline{x}_{i-1}))$ with $k_i\notin \{0,n\}$.	
		
		Let $\underline{x} = [x_n | \ldots | x_1]$ be a redundant cell of height $h$, and let $\underline{z} \neq \underline{x}$ be redundant with $\underline{x} \vdash  \underline{z}$, and $\underline{z}$ is of the form $d_i(\mu(\underline{x}))$ with $i\notin \{0,h+1, n\}$. We will prove that now $\underline{z}$ has height at least $h+1$. For this, let $\underline{y} = \mu(\underline{x}) = [x_n | \ldots | x_{h+2} | a | d | x_h | \ldots | x_1]$ and consider the boundaries $d_i\underline{y}$ for $i \neq h+1$. We distinguish several cases.
		\begin{enumerate}
			\item	$n-1 \geq i \geq h+3$: We have $d_i(\underline{y}) = [x_n | \ldots | x_ix_{i-1} | \ldots | x_{h+2} | a | d | x_h | \ldots |x_1]$, which has height $h+1$ since, as above, $x_k\ldots x_1\in \mathcal{E}$ for $1\leq k \leq h$, $dx_h \ldots x_1 \in \mathcal{E}$ and $adx_h \ldots x_1\notin \mathcal{E}$. As computed above, $\gamma(a, dx_h \ldots x_1)=1$, so $d_i(\underline{y})$ is collapsible.
			
			\item	$h \geq i \geq 1$: For $i \leq h-1$ we have 
			\begin{eqnarray*}
			d_i(\underline{y}) = [x_n | \ldots | x_{h+2} | a | d | x_h | \ldots | x_{i+1}x_{i} | \ldots | x_1], 
			\end{eqnarray*}
and for $i=h$ we have $d_i(\underline{y}) = [x_n | \ldots | x_{h+2} | a | dx_h | x_{h-1} | \ldots | x_1]$. In both cases $d_i(\underline{y})$ has height $h$, because the product of the first $k\leq h$ entries from the right is $x_m\ldots x_1\in \mathcal{E}$ for some $1\leq m\leq h$, or $dx_h \ldots x_1 \in \mathcal{E}$, whereas the product of the first $h+1$ entries from the right gives $adx_h \ldots x_1 \notin \mathcal{E}$. Computing $\gamma(a, dx_h \ldots x_1) = 1$, we see that $d_i(\underline{y})$ is again collapsible.
			
			\item	$i=h+2$: Here, $d_i(\underline{y}) = [ x_n | \ldots | x_{h+3} | x_{h+2}a | d | x_h | \ldots | x_1]$. This cell has height at least $h+1$, for $dx_h \ldots x_1\in \mathcal{E}$ and $x_k\ldots x_1\in\mathcal{E}$ for all $1\leq k\leq h$. The cell $d_i(\underline{y})$ may or may not be redundant.
		\end{enumerate}
		Altogether we have shown that if $\underline{z} \neq \underline{x}$ and $\underline{x} \vdash \underline{z}$, then $\underline{z}$ has strictly larger height than $\underline{x}$. Note that the height of a cell is bounded by its dimension. It follows that every chain $\underline{x}_1 \vdash \underline{x}_2 \vdash \ldots$ eventually stabilizes since all $\underline{x}_n$ have the same dimension.
		
		This finishes the proof of the proposition.
	\end{pf}

\begin{Cor} 
Let $M=M(S)$ be an Artin monoid, $\mathcal{E}$ a generating set closed with respect to left least common multiple and left complement. Then the subcomplex of $BM$ with cells given by $\mathcal{E}_*$ is homotopy equivalent to $BM$. In particular, there is a $\mathbb{Z}$-module complex computing the homology of $M$, with basis $\mathcal{E}_*$ as defined above and differentials given by restriction of the bar differential. 
\end{Cor}

\begin{pf}	
	First, we observe that if $\underline{x} \in \mathcal{E}_*$, then $d_i(\underline{x})$ lies again in $\mathcal{E}_*$ for all $0\leq i\leq n$. Indeed, this is clear for $1\leq i\leq n$. To see this for $d_0[x_n|\ldots | x_1]$, observe that $x_1\in \mathcal{E}$ and $x_k\ldots x_1\in \mathcal{E}$, and their least common multiple is $x_k\ldots x_1$. So $x_k\ldots x_2$ lies in $\mathcal{E}$ as a left complement of elements in $\mathcal{E}$. The corollary now follows from Corollary \ref{UnterkomplexDMT} and Lemma \ref{DMTtoptoalg}. 
\end{pf}

\begin{example}
So far, we were not concerned with the possible choice of a generating system for the results of this section. One important example is provided by the work of Michel (\cite{MichelNoteBraid}, cf. also \cite{BrieskornSaito}): Consider the set of all square-free elements in the monoid $M(S)$, i.e., of all elements which do not have a representative in the free monoid $S^*$ containing a (connected) subword of the form $a^2$. As shown by Tits \cite{Tits}, the square-free elements in an Artin monoid are in one-to-one correspondence with the elements of the associated Coxeter group $W(S)$. In \cite{MichelNoteBraid}, it is shown that the set of square-free elements is closed under left and right least common multiples, as well as under left and right complements. According to the main result of this section, we can thus restrict to the subcomplex of the bar complex generated by $[x_n|\ldots |x_1]$ with $x_n\ldots x_1$ square-free to compute the monoid homology of the Artin monoid.

\end{example}

\bibliographystyle{plain}
\bibliography{Garside}
\end{document}